\documentstyle[amscd, leqno]{amsart}
\theoremstyle{plain}
\newtheorem{theorem}{Theorem}[section]

\newtheorem{lemma}[theorem]{Lemma}

\newtheorem{corollary}[theorem]{Corollary}
\theoremstyle{definition}

\theoremstyle{remark}
\newtheorem{remark}[theorem]{Remark}
\makeatletter

\@addtoreset{equation}{section}
\makeatother

\begin{document}
\title{On the singularity of Quillen metrics}
\author{Ken-Ichi Yoshikawa}
\address{
Graduate School of Mathematical Sciences,
University of Tokyo,
3-8-1 Komaba, Tokyo 153-8914, JAPAN}
\email{yosikawa@@ms.u-tokyo.ac.jp}
\thanks{The author is partially supported 
by the Grants-in-Aid for Scientific Research 
for Encouragement of Young Scientists (B) 16740030, JSPS}

\begin{abstract}
Let $\pi\colon X\to S$ be a holomorphic map
from a compact K\"ahler manifold $(X,g_{X})$ to
a compact Riemann surface $S$.
Let $\Sigma_{\pi}$ be the critical locus of $\pi$
and let $\Delta=\pi(\Sigma_{\pi})$ be 
the discriminant locus.
Let $(\xi,h_{\xi})$ be a holomorphic Hermitian
vector bundle on $X$.
We determine the singularity of the Quillen metric
on $\det R\pi_{*}\xi$ near $\Delta$ 
with respect to $g_{X}|_{TX/S}$ and $h_{\xi}$. 
\end{abstract}

\maketitle


\section
{\bf Introduction}
\par
Let $X$ be a compact K\"ahler manifold 
of dimension $n+1$ with K\"ahler metric $g_{X}$, 
and let $S$ be a compact Riemann surface.
Let $\pi\colon X\to S$ be a surjective 
holomorphic map such that every connected 
component of $X$ is mapped surjectively to $S$.
Let $\Sigma_{\pi}:=\{x\in X;\,d\pi(x)=0\}$
be the critical locus of $\pi$.
For $t\in S$, set $X_{t}:=\pi^{-1}(t)$.
The relative tangent bundle of 
$\pi\colon X\to S$ is the subbundle of 
$TX|_{X\setminus\Sigma_{\pi}}$ defined as
$TX/S:=\ker\pi_{*}|_{X\setminus\Sigma_{\pi}}$.
Set 
$$
\Delta:=\pi(\Sigma_{\pi}),
\qquad
S^{o}:=S\setminus\Delta,
\qquad
X^{o}:=X|_{S^{o}},
\qquad
\pi^{o}:=\pi|_{X^{o}}.
$$
Then $\pi^{o}\colon X^{o}\to S^{o}$ is 
a holomorphic family of compact K\"ahler 
manifolds.
Let $g_{X/S}:=g_{X}|_{TX/S}$ be the 
Hermitian metric on $TX/S$ induced from 
$g_{X}$. 
\par
Let $\xi\to X$ be a holomorphic vector bundle 
on $X$ equipped with 
a Hermitian metric $h_{\xi}$.
Let $\lambda(\xi)=\det R\pi_{*}\xi$ be 
the determinant of the cohomologies of $\xi$.
By \cite{BGS88}, \cite{Quillen85}, 
\cite{Soule92}, 
$\lambda(\xi)|_{S^{o}}$ is equipped 
with the Quillen metric 
$\|\cdot\|^{2}_{\lambda(\xi),Q}$
with respect to the metrics
$g_{X/S}$ and $h_{\xi}$. 
\par
Let $0\in\Delta$ be an arbitrary critical 
value of $\pi$,
and let $({\cal U},t)$ be a coordinate 
neighborhood of $S$ centered at $0$
with ${\cal U}\cap\Delta=\{0\}$.
Set ${\cal U}^{o}:={\cal U}\setminus\{0\}$.
\par
Let $\sigma$ be a nowhere vanishing 
holomorphic section of $\lambda(\xi)$ 
on $\cal U$.
Then $\log\|\sigma\|^{2}_{\lambda(\xi),Q}$ 
is a $C^{\infty}$ function on 
${\cal U}^{o}$ by \cite{BGS88}.
The purpose of this article is to study 
the behavior of 
$\log\|\sigma(t)\|^{2}_{\lambda(\xi),Q}$ 
as $t\to0$.
\par
For a holomorphic vector bundle $F$ over 
a complex manifold with zero-section $Z$, 
define the projective-space bundle 
${\Bbb P}(F)$ as 
${\Bbb P}(F):=(F\setminus Z)/{\Bbb C}^{*}$.
The dual projective-space bundle 
${\Bbb P}(F)^{\lor}$ is defined as 
${\Bbb P}(F)^{\lor}:={\Bbb P}(F^{\lor})$,
where $F^{\lor}$ is the dual vector bundle
of $F$.
\par
Following Bismut \cite{Bismut97}, 
we consider the Gauss map 
$\mu\colon X\setminus\Sigma_{\pi}\to
{\Bbb P}(TX)^{\lor}$ that assigns
$x\in X\setminus\Sigma_{\pi}$ the hyperplane
$\ker(\pi_{*})_{x}\in{\Bbb P}(T_{x}X)^{\lor}$.
Since $\mu$ extends to a meromorphic map
$\mu\colon X\dashrightarrow{\Bbb P}(TX)^{\lor}$,
there exists a resolution 
$q\colon (\widetilde{X},E)\to(X,\Sigma_{\pi})$
of the indeterminacy of $\mu$ such that
$\widetilde{\mu}:=\mu\circ q$ extends to 
a holomorphic map from $\widetilde{X}$ to 
${\Bbb P}(TX)^{\lor}$
and such that $E$ is a normal crossing divisor 
of $\widetilde{X}$.
(For the scheme structure of $E$, 
see Sect.\,3.)
Let $U$ be the universal hyperplane bundle of 
rank $n=\dim X/S$ over ${\Bbb P}(TX)^{\lor}$, 
and let $H:={\cal O}_{{\Bbb P}(TX)^{\lor}}(1)$.
\par
After Barlet \cite{Barlet82}, we define
a subspace of $C^{0}({\cal U})$ by
$$
{\cal B}({\cal U})
:=
C^{\infty}({\cal U})
\oplus
\bigoplus_{r\in{\Bbb Q}\cap(0,1]}
\bigoplus_{k=0}^{n}
|t|^{2r}(\log|t|)^{k}\cdot C^{\infty}({\cal U}).
$$
A function
$\varphi(t)\in{\cal B}({\cal U})$ has 
an asymptotic expansion at $0\in\Delta$,
i.e., there exist
$r_{1},\ldots,r_{m}\in{\Bbb Q}\cap(0,1]$ and
$f_{0},f_{l,k}\in C^{\infty}({\cal U})$,
$l=1,\ldots,m$, $k=0,\ldots,n$, such that
$$
\varphi(t)
=
f_{0}(t)+
\sum_{l=1}^{m}\sum_{k=0}^{n}
|t|^{2r_{l}}(\log|t|)^{k}\,f_{l,k}(t).
$$
\par
In what follows,
if $f(t),g(t)\in 
C^{\infty}({\cal U}^{o})$ satisfies
$f(t)-g(t)\in{\cal B}({\cal U})$, 
we write 
$$
f\equiv_{\cal B}g.
$$ 
\par
For a complex vector bundle $F$ over
a complex manifold, $c_{i}(F)$,
${\rm Td}(F)$, and ${\rm ch}(F)$ denote
the $i$-th Chern class, the Todd genus,
and the Chern character of $F$, respectively.
\par
We can state the main result of 
this article, which generalizes 
\cite[\S5]{Bismut97} and \cite{Yoshikawa98}:

\begin{theorem}
The following identity holds:
$$
\log\|\sigma\|^{2}_{Q,\lambda(\xi)}
\equiv_{\cal B}
\left(\int_{E\cap q^{-1}(X_{0})}
\widetilde{\mu}^{*}
\left\{{\rm Td}(U)\,
\frac{{\rm Td}(H)-1}{c_{1}(H)}\right\}\,
q^{*}{\rm ch}(\xi)\right)\log|t|^{2}.
$$
\end{theorem}

By Theorem 1.1, $\|\cdot\|^{2}_{Q,\lambda(\xi)}$ 
extends to a singular Hermitian metric on 
$\lambda(\xi)$. Let $\pi_{*}$ denote 
the integration along the fibers of $\pi$.
As a consequence of Theorem 1.1
and the curvature formula for Quillen metrics
\cite{BGS88}, we get the following:

\begin{corollary}
The $(1,1)$-form 
$\pi_{*}({\rm Td}(TX/S,g_{X/S})\,
{\rm ch}(\xi,h_{\xi}))^{(1,1)}$ lies in
$L^{p}_{\rm loc}(S)$ for some $p>1$,
and the curvature current of
$(\lambda(\xi),\|\cdot\|_{Q,\lambda(\xi)})$
is given by the following formula on ${\cal U}$:
$$
\begin{aligned}
c_{1}(\lambda(\xi),\|\cdot\|_{Q,\lambda(\xi)})
&
=
\pi_{*}({\rm Td}(TX/S,g_{X/S})\,
{\rm ch}(\xi,h_{\xi}))^{(1,1)}
\\
&\quad-
\left(\int_{E\cap q^{-1}(X_{0})}
\widetilde{\mu}^{*}
\left\{{\rm Td}(U)\,
\frac{{\rm Td}(H)-1}{c_{1}(H)}\right\}\,
q^{*}{\rm ch}(\xi)\right)\,\delta_{0},
\end{aligned}
$$
where $\delta_{0}$ denotes the Dirac 
$\delta$-current supported at $0$.
\end{corollary}

The proof of Theorem 1.1 is quite similar
to that of Bismut in \cite[\S5]{Bismut97},
and we just follow his argument. 
There are essentially no new ideas except
a systematic use of the Gauss maps
for the family $\pi\colon X\to S$;
in fact, the Gauss maps were already
used by Bismut in \cite{Bismut97}.
\par
The existence of an asymptotic expansion
of the Quillen norm 
$\log\|\sigma\|^{2}_{Q,\lambda(\xi)}$
was first shown
by Bismut-Bost\cite[Sect.\,13.(b)]{BB90}
when $\pi\colon X\to S$ is a family of curves
and by the author \cite{Yoshikawa98}
when $\Sigma_{\pi}$ is isolated.
In \cite{FangLuYoshikawa05}, 
Theorem 1.1 shall play an crucial role
in the study of analytic torsion of Calabi-Yau
threefolds.
\par
Let ${\bf s}_{\Delta}$ be a section of 
${\cal O}_{S}(\Delta)$ defining the reduced 
divisor $\Delta$.
Let $\|\cdot\|$ be a $C^{\infty}$ 
Hermitian metric on ${\cal O}_{S}(\Delta)$.
By Theorem 1.1, 
$$
\log\|\sigma(t)\|^{2}_{Q,\lambda(\xi)}
-
\left(\int_{E\cap q^{-1}(X_{0})}
\widetilde{\mu}^{*}
\left\{{\rm Td}(U)\,
\frac{{\rm Td}(H)-1}{c_{1}(H)}\right\}\,
q^{*}{\rm ch}(\xi)\right)
\log\|{\bf s}_{\Delta}(t)\|^{2}
$$
has a finite limit as $t\to0$. In Section 6,
we shall compute this limit in terms of 
various secondary objects,
which extends some results in
\cite[\S5]{Bismut97}.
\par
This article is organized as follows.
In Sections 2 and 3,
we explain the Gauss maps
associated to the family $\pi\colon X\to S$
and their resolutions.
In Sections 5 and 6,
we prove the main theorem.
In Sections 7 and 8, 
we verify the compatibility of
Theorem 1.1 with the corresponding 
earlier results of 
Bismut \cite{Bismut97}
and the author \cite{Yoshikawa98}.
In Sections 4 and 9, 
we prove some technical results.
The problem treated in Section 9
seems to be related with 
the regularity problem of 
the star products of Green currents
\cite{BostGilletSoule94}.
\par
For a complex manifold, we set
$d^{c}=
\frac{1}{4\pi i}(\partial-\bar{\partial})$.
Hence
$dd^{c}=
\frac{1}{2\pi i}\bar{\partial}\partial$.
We keep the notation in Sect.\,1 throughout 
this article.


\section
{\bf The Gauss maps}
\par
Let $\Omega^{1}_{X}$ be the holomorphic 
cotangent bundle of $X$.
Let 
$\varPi\colon
{\Bbb P}(\Omega^{1}_{X}\otimes\pi^{*}TS)\to X$ 
be the projective-space bundle associated with 
$\Omega^{1}_{X}\otimes\pi^{*}TS$. 
Since $\dim S=1$, we have 
${\Bbb P}(\Omega^{1}_{X}\otimes\pi^{*}TS)
=
{\Bbb P}(\Omega^{1}_{X})$.
Let 
$\varPi^{\lor}\colon{\Bbb P}(TX)^{\lor}\to X$
be the dual projective-space bundle of 
${\Bbb P}(TX)$,
whose fiber ${\Bbb P}(T_{x}X)^{\lor}$ 
is the set of hyperplanes of $T_{x}X$ 
passing through the zero vector of $T_{x}X$. 
We have the canonical isomorphisms
$$
{\Bbb P}(\Omega^{1}_{X}\otimes\pi^{*}TS)
=
{\Bbb P}(\Omega^{1}_{X})
\cong
{\Bbb P}(TX)^{\lor}.
$$
\par
Let $x\in X\setminus\Sigma_{\pi}$.
Let $t$ be a holomorphic local coordinate 
of $S$ near $\pi(x)\in S$. 
We define the Gauss maps
$\nu\colon X\setminus\Sigma_{\pi}\to
{\Bbb P}(\Omega^{1}_{X}\otimes\pi^{*}TS)$ 
and
$\mu\colon 
X\setminus\Sigma_{\pi}\to{\Bbb P}(TX)^{\lor}$
by
$$
\nu(x)
:=
[d\pi_{x}]
=
\left[
\sum_{i=0}^{n}\frac{\partial(t\circ\pi)}
{\partial z_{i}}(x)\,dz_{i}\otimes
\frac{\partial}{\partial t}
\right],
\qquad
\mu(x)
:=
[T_{x}X_{\pi(x)}].
$$
Under the canonical isomorphism
${\Bbb P}(\Omega^{1}_{X}\otimes\pi^{*}TS)
\cong
{\Bbb P}(TX)^{\lor}$, 
one has
$$
\nu=\mu.
$$
\par
Let 
$$
L
:=
{\cal O}_{{\Bbb P}
(\Omega^{1}_{X}\otimes\pi^{*}TS)}
(-1)\subset
\varPi^{*}(\Omega^{1}_{X}\otimes\pi^{*}TS)
$$
be the tautological line bundle over 
${\Bbb P}(\Omega^{1}_{X}\otimes\pi^{*}TS)$, 
and set
$$
Q
:=
\varPi^{*}(\Omega^{1}_{X}\otimes\pi^{*}TS)/L.
$$ 
We have the exact sequence
of holomorphic vector bundles 
on ${\Bbb P}(\Omega^{1}_{X}\otimes\pi^{*}TS)$:
$$
{\cal S}\colon
0
\longrightarrow
L
\longrightarrow 
\varPi^{*}(\Omega^{1}_{X}\otimes\pi^{*}TS)
\longrightarrow
Q
\longrightarrow
0.
$$
\par
Let
$H={\cal O}_{{\Bbb P}(T_{X})^{\lor}}(1)$,
and let $U$ be the universal hyperplane 
bundle of $(\varPi^{\lor})^{*}T{X}$. 
Then the dual of ${\cal S}$ is given by
$$
{\cal S}^{\lor}\colon
0
\longrightarrow 
U
\longrightarrow
(\varPi^{\lor})^{*}TX
\longrightarrow
H
\longrightarrow 
0.
$$
Since 
$T_{x}X_{\pi(x)}=\{v\in T_{x}X;\,
d\pi_{x}(v)=0\}$, 
we have on $X\setminus\Sigma_{\pi}$
$$
TX/S=\mu^{*}U.
$$
\par
Let $g_{U}$ be the Hermitian metric on $U$ 
induced from $(\varPi^{\lor})^{*}g_{X}$, 
and 
let $g_{H}$ be the Hermitian metric on $H$ 
induced from $(\varPi^{\lor})^{*}g_{X}$ 
by the $C^{\infty}$-isomorphism 
$H\cong U^{\perp}$.
On $X\setminus\Sigma_{\pi}$, we have
$$
(TX/S,g_{{\cal X}/S})
=
\mu^{*}(U,g_{U}).
$$
\par
Let $g_{S}$ be a Hermitian metric on $S$.
Let $g_{\Omega^{1}_{X}}$ be the Hermitian 
metric on $\Omega^{1}_{X}$ 
induced from $g_{X}$.
Let $g_{L}$ be the Hermitian metric on $L$ 
induced from the metric
$\varPi^{*}(g_{\Omega^{1}_{X}}\otimes
\pi^{*}g_{S})$ 
by the inclusion 
$L\subset\varPi^{*}
(\Omega^{1}_{X}\otimes\pi^{*}TS)$. 
Let $g_{Q}$ be the Hermitian metric on $Q$ 
induced from
$\varPi^{*}(g_{\Omega^{1}_{X}}\otimes
\pi^{*}g_{S})$ 
by the $C^{\infty}$-isomorphism 
$Q\cong L^{\perp}$.
\par
Let $c_{1}(L,g_{L})$ be the Chern form of $(L,g_{L})$.
Since $d\pi$ is a nowhere vanishing 
holomorphic section of 
$\nu^{*}L|_{{X}\setminus\Sigma_{\pi}}$, 
we get the following equation on 
$X\setminus\Sigma_{\pi}$
$$
-dd^{c}\log\|d\pi\|^{2}
=
\nu^{*}c_{1}(L,g_{L}).
$$


\section
{\bf Resolution of the Gauss maps}
\par
Since $\Sigma_{\pi}$ is a proper analytic 
subset of $X$,
the maps 
$\nu\colon X\setminus\Sigma_{\pi}\to
{\Bbb P}(\Omega^{1}_{X}\otimes\pi^{*}TS)$ 
and
$\mu\colon X\setminus\Sigma_{\pi}\to
{\Bbb P}(TX)^{\lor}$ 
extend to meromorphic maps
$\nu\colon X\dashrightarrow
{\Bbb P}(\Omega^{1}_{X}\otimes\pi^{*}TS)$
and
$\mu\colon X\dashrightarrow{\Bbb P}(TX)^{\lor}$
by \cite[Th.\,4.5.3]{NoguchiOchiai90}.
By Hironaka, there exists 
a compact K\"ahler manifold $\widetilde{X}$, 
a normal crossing divisor 
$E\subset\widetilde{X}$, 
a birational holomorphic map
$q\colon\widetilde{X}\to{X}$, 
and holomorphic maps
$\widetilde{\nu}\colon\widetilde{X}\to
{\Bbb P}(\Omega^{1}_{X}\otimes\pi^{*}TS)$ 
and
$\widetilde{\mu}\colon\widetilde{X}\to
{\Bbb P}(TX)^{\lor}$ 
satisfying the following conditions:
\newline{(i)}
$q
|_{\widetilde{X}\setminus q^{-1}(\Sigma_{\pi})}
\colon
\widetilde{X}\setminus q^{-1}(\Sigma_{\pi})
\to X\setminus\Sigma_{\pi}$ 
is an isomorphism;
\newline{(ii)}
$q^{-1}(\Sigma_{\pi})=E$;
\newline{(iii)}
$(\pi\circ q)^{-1}(b)$ 
is a normal crossing divisor
of $\widetilde{X}$ for all $b\in\Delta$;
\newline{(iv)}
$\widetilde{\nu}=\nu\circ q$ 
and
$\widetilde{\mu}=\mu\circ q$ 
on $\widetilde{X}\setminus E$.
\newline
Then $\widetilde{\nu}=\widetilde{\mu}$
under the canonical isomorphism
${\Bbb P}(\Omega^{1}_{X}\otimes\pi^{*}TS)
\cong{\Bbb P}(TX)^{\lor}$.
We set 
$$
\widetilde{\pi}:=\pi\circ q
$$ 
and
$\widetilde{X}_{s}:=\widetilde{\pi}^{-1}(s)$
for $s\in S$. Similarly, we set
$E_{b}:=E\cap\widetilde{X}_{b}$ for $b\in\Delta$.
Sinec $E=q^{-1}(\Sigma_{\pi})\subset
\widetilde{\pi}^{-1}(\Delta)$, we have
$E=\amalg_{b\in\Delta}E_{b}$.
\par
Let ${\cal I}_{\Sigma_{\pi}}$ be 
the ideal sheaf of $\Sigma_{\pi}$. 
For every $p\in\Sigma_{\pi}$,
the sheaf ${\cal I}_{\Sigma_{\pi}}$ 
has the following expression
on a neighborhood of $p$:
$$
{\cal I}_{\Sigma_{\pi}}
=
{\cal O}_{X}
\left(
\frac{\partial(t\circ\pi)}{\partial z_{0}}(z)
,\cdots,
\frac{\partial(t\circ\pi)}{\partial z_{n}}(z)
\right).
$$
Define the ideal sheaf
${\cal I}_{E}$ of $E$ as
$$
{\cal I}_{E}=q^{-1}{\cal I}_{\Sigma_{\pi}}.
$$
\par
Denote by $\delta_{E}$ the $(1,1)$-current on
$\widetilde{X}$ defined as the integration 
over $E$, i.e.,
$\delta_{E}(\psi):=\int_{E}\psi|_{E}$ 
for all $C^{\infty}$ $(n,n)$-form 
on $\widetilde{X}$.
Since $\widetilde{\nu}^{*}L=q^{*}\nu^{*}L$,
$q^{*}d\pi$ extends to a holomorphic section 
of $\widetilde{\nu}^{*}L$ with 
zero divisor $E$ by the definition
of the ideal sheaf ${\cal I}_{E}$. 
By the Poincar\'e-Lelong formula, 
the following identity of currents on 
$\widetilde{X}$ holds
$$
-dd^{c}(q^{*}\log\|d\pi\|^{2})
=
\widetilde{\nu}^{*}c_{1}(L,g_{L})-\delta_{E}.
$$


\section
{\bf Regularity of the direct image of 
differential forms}
\par
Recall that $({\cal U},t)$ is
a coordinate neighborhood of $S$ centered
at the critical value $0\in\Delta$.
Set 
$D:=\{(s,t)\in S\times{\cal U};\,s=t\}$.
Then $D$ is a divisor of $S\times{\cal U}$.
Let $[D]$ be the line bundle on 
$S\times{\cal U}$ 
defined by the divisor $D$.
Let ${\bf s}_{D}$ be a section of $[D]$ 
with zero divisor $D$.
Let $B\subset S$ be a finite subset with
$0\in B$. By shrinking ${\cal U}$ if necessary,
we may assume that ${\cal U}\cap B=\{0\}$.
Let $\|\cdot\|_{D}$ be a $C^{\infty}$ 
Hermitian metric on $[D]$ such that
\begin{equation}
\|{\bf s}_{D}(b,t)\|_{D}=1,
\qquad
\forall\,(b,t)\in(B\setminus\{0\})\times{\cal U}.
\end{equation}
We set
${\bf s}_{t}:={\bf s}_{D}|_{S\times\{t\}}$ 
and
$\|\cdot\|_{t}:=\|\cdot\|_{D}|_{S\times\{t\}}$
for $t\in{\cal U}$. 
Then ${\rm div}({\bf s}_{t})=\{t\}$ and
$\|{\bf s}_{t}\|_{t}^{2}\in 
C^{\infty}(S\times{\cal U})$.
\par
Let $V$ be a compact connected complex manifold
with $\dim V=n+1$.
Let $f\colon V\to S$ be a proper surjective 
holomorphic map. We set $V_{t}:=f^{-1}(t)$
for $t\in S$.
\par
Let $\overline{F}:=(F,\|\cdot\|_{F})$ be 
a holomorphic Hermitian line bundle on $V$,
and let $\alpha$ be a holomorphic section of $F$ with
$$
{\rm div}(\alpha)\subset\sum_{b\in B}V_{b}.
$$
\par
Denote by $f_{*}$ 
the integration along the fibers of $f$. 
In Section 4, we assume that
$\varphi$ is a $\partial$-closed and 
$\bar{\partial}$-closed $C^{\infty}$ 
$(n,n)$-form on $V$.

\begin{lemma}
There exists a H\"order continuous function
$\eta$ on $\cal U$ such that
$$
f_{*}
\{(\log\|\alpha\|_{F}^{2})\,\varphi\}^{(0,0)}
-
\left(
\int_{{\rm div}(\alpha)\cap V_{0}}\varphi
\right)\,
\log\|{\bf s}_{0}\|_{0}^{2}
=
\eta.
$$
\end{lemma}

\begin{pf}
Since $\log\|\alpha\|_{F}^{2}\,\varphi$ is 
a locally integrable differential form on 
$V$, we have
$f_{*}\{(\log\|\alpha\|^{2})\,\varphi\}^{(0,0)}
\in L^{1}_{\rm loc}(S)\cap C^{\infty}(S^{o})$.
Since $dd^{c}$ commutes with $f_{*}$ 
and since $\varphi$ is $d$ and $d^{c}$-closed,
we get the following equation of currents 
on $\cal U$:
\begin{equation}
\begin{aligned}
dd^{c}f_{*}
\{(\log\|\alpha\|_{F}^{2})\,\varphi\}^{(0,0)}
&
=
[f_{*}\{
dd^{c}((\log\|\alpha\|_{F}^{2})\wedge\varphi)
\}]^{(1,1)}
\\
&
=
-[f_{*}\left\{
(c_{1}(\overline{F})-\delta_{{\rm div}(\alpha)})
\wedge\varphi\right\}]^{(1,1)}
\\
&
=
\left(
\int_{{\rm div}(\alpha)\cap V_{0}}\varphi
\right)\,
\delta_{0}-
[f_{*}\{
c_{1}(\overline{F})\wedge\varphi\}]^{(1,1)}.
\end{aligned}
\end{equation}
By Lemma 9.2 below, there exists 
$\psi\in{\cal B}({\cal U})$ such that
$$
[f_{*}\{
c_{1}(\overline{F})
\wedge\varphi\}]^{(1,1)}(t)
=
\psi(t)\,
\frac{dt\wedge d\bar{t}}{|t|^{2}},
\qquad
\psi(0)=0.
$$
Since $\psi(0)=0$, there exists
$\nu\in{\Bbb Q}\cap(0,1]$ 
such that
$\psi(t)\in
\sum_{k\leq n}
|t|^{2\nu}(\log|t|)^{k}\cdot{\cal B}({\cal U})$.
Hence
$|t|^{-2}\psi(t)\in L^{p}_{\rm loc}({\cal U})$
for some $p>1$.
By the ellipticity of the Laplacian
and the Sobolev embedding theorem, 
there exists a H\"older continuous function 
$\chi$ on $\cal U$ satisfying 
the following equation of currents on ${\cal U}$
$$
[f_{*}\{c_{1}(\overline{F})
\wedge\varphi\}]^{(1,1)}
=
dd^{c}\chi.
$$
This, together with (4.2) and the equation
of currents 
$dd^{c}\log|t|^{2}=\delta_{0}$ 
on $\cal U$, implies the assertion, 
because 
$\log\|{\bf s}_{0}\|_{0}^{2}-\log|t|^{2}
\in C^{\infty}({\cal U})$.
\end{pf}

\begin{lemma}
The following identity holds for all
$t\in{\cal U}^{o}$:
$$
\begin{aligned}
\int_{V_{t}}
(\log\|\alpha\|_{F}^{2})\,\varphi
&
=
\left(
\int_{{\rm div}(\alpha)\cap V_{0}}\varphi
\right)\,
\log\|{\bf s}_{t}(0)\|_{t}^{2}
-
\int_{V}
(f^{*}\log\|{\bf s}_{t}\|_{t}^{2})\,
c_{1}(\overline{F})\wedge\varphi
\\
&
\quad
+
\int_{V}
(\log\|\alpha\|_{F}^{2})\,
f^{*}c_{1}([t],\|\cdot\|_{t})
\wedge\varphi.
\end{aligned}
$$
\end{lemma}

\begin{pf}
Since $V_{t}\cap{\rm div}(\alpha)=\emptyset$
for $t\in{\cal U}^{o}$, $V_{t}$ meets ${\rm div}(\alpha)$
properly. Since $\varphi$ is $\partial$ and 
$\bar{\partial}$-closed,
we deduce from \cite[Th.\,2.2.2]{GilletSoule90b}
the following identity by setting
$X=W=V$,
$Y=V_{t}$,
$Z={\rm div}(\alpha)$,
and
$g_{Y}=-f^{*}
\log\|{\bf s}_{t}\|_{t}^{2}$,
$g_{Z}=-\log\|\alpha\|_{F}^{2}$
in
\cite[Sect.\,2.2.2]{GilletSoule90b}:
\begin{equation}
\begin{aligned}
\int_{V_{t}}
(\log\|\alpha\|_{F}^{2})\,\varphi
&
=
\sum_{b\in B}
\left(
\int_{{\rm div}(\alpha)\cap V_{b}}\varphi
\right)\,
\log\|{\bf s}_{t}(b)\|_{t}^{2}
-
\int_{V}
(f^{*}\log\|{\bf s}_{t}\|_{t}^{2})\,
c_{1}(\overline{F})\wedge\varphi
\\
&
\quad
+
\int_{V}
(\log\|\alpha\|_{F}^{2})\,
f^{*}c_{1}([t],\|\cdot\|_{t})
\wedge\varphi,
\end{aligned}
\end{equation}
where we used the assumption
${\rm div}(\alpha)\subset\sum_{b\in B}V_{b}$.
(See also \cite[p.59, l.3-l.7]{Soule92}.)
Since $\|{\bf s}_{t}(b)\|_{t}=1$ for
$(b,t)\in(B\setminus\{0\})\times{\cal U}$ by (4.1),
the result follows from (4.3).
\end{pf}

\begin{lemma}
The following identity holds
$$
\begin{aligned}
\,&
\lim_{t\to0}
\left\{
\int_{V_{t}}
(\log\|\alpha\|_{F}^{2})\,\varphi
-\left(
\int_{{\rm div}(\alpha)\cap V_{0}}\varphi
\right)\,
\log\|{\bf s}_{0}(t)\|_{0}^{2}
\right\}=
\\
&
\int_{V}
(\log\|\alpha\|_{F}^{2})\,
f^{*}c_{1}([0],\|\cdot\|_{0})\wedge
\varphi
-
\int_{V}
(f^{*}\log\|{\bf s}_{0}\|_{0}^{2})\,
c_{1}(\overline{F})\wedge\varphi.
\end{aligned}
$$
\end{lemma}

\begin{pf}
By Lemma 4.2, we have
\begin{equation}
\begin{aligned}
\int_{V_{t}}
(\log\|\alpha\|_{F}^{2})\,\varphi
&
=
\left(
\int_{{\rm div}(\alpha)\cap V_{0}}\varphi
\right)\,
\log\|{\bf s}_{0}(t)\|_{0}^{2}
-
\int_{V}
(f^{*}\log\|{\bf s}_{t}\|_{t}^{2})\,c_{1}(\overline{F})
\wedge\varphi
\\
&
\quad
+
\int_{V}
(\log\|\alpha\|_{F}^{2})\,f^{*}c_{1}([t],\|\cdot\|_{t})
\wedge\varphi
+
\left(
\int_{{\rm div}(\alpha)\cap V_{0}}\varphi
\right)\,
\log\frac{\|{\bf s}_{t}(0)\|_{t}^{2}}
{\|{\bf s}_{0}(t)\|_{0}^{2}}.
\end{aligned}
\end{equation}
Since 
$\lim_{s\to0}\log(\|{\bf s}_{t}(0)\|_{t}^{2}/
\|{\bf s}_{0}(t)\|_{0}^{2})=0$, 
the assertion follows from (4.4).
\end{pf}

\begin{lemma}
The following identity of functions 
on ${\cal U}^{o}$ hold:
$$
f_{*}
\{(\log\|\alpha\|_{F}^{2})\,\varphi\}^{(0,0)}
\equiv_{\cal B}
\left(
\int_{{\rm div}(\alpha)\cap V_{0}}\varphi
\right)\,
\log\|{\bf s}_{0}\|_{0}^{2}.
$$
\end{lemma}

\begin{pf}
For $t\in{\cal U}^{o}$, set
$$
I_{1}(t)
:=
\int_{V}
(f^{*}\log\|{\bf s}_{t}\|_{t}^{2})\,
c_{1}(\overline{F})\,\varphi,
\qquad
I_{2}(t)
:=
\int_{V}
(\log\|\alpha\|_{F}^{2})\,
f^{*}c_{1}([t],\|\cdot\|_{t})\,
\varphi.
$$
By (4.4), it suffices to prove that
$I_{1}\in{\cal B}({\cal U})$ and
$I_{2}\in{\cal B}({\cal U})$.
\par
Let 
$\{(W_{\lambda},z_{\lambda})\}_{\lambda\in\Lambda}$ 
be a system of local coordinates on $V$.
Since $V$ is compact, 
we may assume $\#\Lambda<+\infty$.
For every $\lambda\in\Lambda$, there exist
$F_{\lambda}\in{\cal O}(W_{\lambda})$,
$G_{\lambda}\in{\cal O}(W_{\lambda})$, 
$A_{\lambda}\in C^{\infty}(W_{\lambda})$,
and
$B_{\lambda}\in 
C^{\infty}(W_{\lambda}\times{\cal U})$ 
such that
$$
\widetilde{\pi}^{*}\log\|{\bf s}_{t}\|_{t}^{2}
|_{W_{\lambda}}(z_{\lambda})
=
\log|F_{\lambda}(z_{\lambda})-t|^{2}
+
B_{\lambda}(z_{\lambda},t),
$$
$$
\log\|\alpha\|_{F}^{2}
|_{W_{\lambda}}(z_{\lambda})
=
\log|G_{\lambda}(z_{\lambda})|^{2}
+
A_{\lambda}(z_{\lambda}).
$$
\par
Let 
$\{\varrho_{\lambda}\}_{\lambda\in\Lambda}$ 
be a partition of unity of $V$ 
subject to the covering
$\{W_{\lambda}\}_{\lambda\in\Lambda}$. 
We set 
$\chi_{\lambda}
:=
\varrho_{\lambda}\,c_{1}(\overline{F})\,\varphi$.
Then
\begin{equation}
I_{1}(t)
=
\sum_{\lambda\in\Lambda}\int_{W_{\lambda}}
\log|F_{\lambda}(z_{\lambda})-t|^{2}
\cdot\chi_{\lambda}(z_{\lambda})+
\sum_{\lambda\in\Lambda}\int_{W_{\lambda}}
B_{\lambda}(z_{\lambda},t)
\,\chi_{\lambda}(z_{\lambda}).
\end{equation}
Since the first term of the right hand side 
of (4.5) lies in ${\cal B}({\cal U})$ 
by Theorem 9.1 below,
we get $I_{1}\in{\cal B}({\cal U})$.
\par
We set $\theta_{\lambda}
:=
\varrho_{\lambda}\,
\widetilde{\pi}^{*}c_{1}([t],\|\cdot\|_{t})\,
\varphi$. 
Then 
$\theta_{\lambda}(z_{\lambda},t)$ 
is a $C^{\infty}$ $(n+1,n+1)$-form 
on $W_{\lambda}\times{\cal U}$. Since
$$
I_{2}(t)
=
\sum_{\lambda\in\Lambda}\int_{W_{\lambda}}
\log|G_{\lambda}(z_{\lambda})|^{2}
\cdot\theta_{\lambda}(z_{\lambda},t)+
\sum_{\lambda\in\Lambda}\int_{W_{\lambda}}
A_{\lambda}(z_{\lambda})\,
\theta_{\lambda}(z_{\lambda},t),
$$ 
we get $I_{2}\in C^{\infty}({\cal U})$.
This completes the proof.
\end{pf}

\begin{corollary}
The following identity holds
$$
\begin{aligned}
\,&
\lim_{t\to0}
\left\{
\int_{\widetilde{X}_{t}}
q^{*}(\log\|d\pi\|^{2})\,\varphi
-\left(\int_{E_{0}}\varphi\right)\,
\log\|{\bf s}_{0}(t)\|_{0}^{2}
\right\}=
\\
&
\int_{\widetilde{X}}
(q^{*}\log\|d\pi\|^{2})\,
\widetilde{\pi}^{*}
c_{1}([0],\|\cdot\|_{0})\wedge
\varphi
-
\int_{\widetilde{X}}
(\widetilde{\pi}^{*}
\log\|{\bf s}_{0}\|_{0}^{2})\,
\widetilde{\nu}^{*}c_{1}(L,g_{L})\wedge\varphi.
\end{aligned}
$$
\end{corollary}

\begin{pf}
Setting $V=\widetilde{X}$, $f=\widetilde{\pi}$,
$\overline{F}=\widetilde{\nu}^{*}(L,g_{L})$ and
$\alpha=q^{*}(d\pi)$ in Lemma 4.3,
we get the result.
\end{pf}

\begin{corollary}
The following identity of functions 
on ${\cal U}^{o}$ hold:
$$
\widetilde{\pi}_{*}
(q^{*}(\log\|d\pi\|^{2})\,\varphi)^{(0,0)}
\equiv_{\cal B}
\left(\int_{E_{0}}\varphi\right)\,
\log\|{\bf s}_{0}\|_{0}^{2}.
$$
\end{corollary}

\begin{pf}
Setting $V=\widetilde{X}$, $f=\widetilde{\pi}$,
$\overline{F}=\widetilde{\nu}^{*}(L,g_{L})$ and
$\alpha=q^{*}(d\pi)$ in Lemma 4.4,
we get the result.
\end{pf}


\section
{\bf Behavior of the Quillen norm of 
the Knudsen-Mumford section}
\par
Let $\Gamma\subset X\times S$ be 
the graph of $\pi$, which is a smooth divisor 
on $X\times S$. 
Let $[\Gamma]$ be the holomorphic line bundle 
on $X\times S$ associated to $\Gamma$. 
Let $s_{\Gamma}\in H^{0}(X\times S,[\Gamma])$ 
be the canonical section of $[\Gamma]$,
so that ${\rm div}(s_{\Gamma})=\Gamma$.
We identify $X$ with $\Gamma$.
\par
Let 
$i\colon\Gamma\hookrightarrow X\times S$
be the inclusion.
Let 
$p_{1}\colon X\times S\to X$ 
and
$p_{2}\colon X\times S\to S$ 
be the projections.
On $X\times S$, 
we have the exact sequence of
coherent sheaves,
\begin{equation}
0
\longrightarrow
{\cal O}_{X\times S}
([\Gamma]^{-1}\otimes p_{1}^{*}\xi)
@>\otimes s_{\Gamma}>>
{\cal O}_{X\times S}(p_{1}^{*}\xi)
\longrightarrow
i_{*}{\cal O}_{\Gamma}(p_{1}^{*}\xi)
\longrightarrow0.
\end{equation}
\par
Let 
$\lambda(p_{1}^{*}\xi)$, 
$\lambda
([\Gamma]^{-1}\otimes p_{1}^{*}\xi)$, 
$\lambda(\xi)$ 
be the determinants of the direct images 
$R(p_{2})_{*}p_{1}^{*}\xi$, 
$R(p_{2})_{*}
([\Gamma]^{-1}\otimes p_{1}^{*}\xi)$, 
$R\pi_{*}\xi$, 
respectively. 
By definition
\cite{BGS88},
\cite{KnudsenMumford76},
\cite{Soule92},
$$
\lambda(\xi)
=
\bigotimes_{q\geq0}
(\det R^{q}\pi_{*}\xi)^{(-1)^{q}}.
$$
Under the isomorphism
$p_{1}^{*}\xi|_{\Gamma}\cong\xi$ 
induced from the identification 
$p_{1}\colon\Gamma\to X$, 
the holomorphic line bundle on $S$
$$
\lambda
:=
\lambda
\left([\Gamma]^{-1}\otimes p_{1}^{*}\xi\right)
\otimes
\lambda(p_{1}^{*}\xi)^{-1}
\otimes
\lambda(\xi)
$$
carries the canonical nowhere vanishing
holomorphic section $\sigma_{KM}$
by \cite{BismutLebeau91}, \cite{KnudsenMumford76}.
\par
Let ${\cal V}\subset{\cal U}$ be 
a relatively compact neighborhood 
of $0\in\Delta$, 
and set
${\cal V}^{o}:={\cal V}\setminus\{0\}$. 
On $\pi^{-1}({\cal U})$, 
we identify $\pi$ (resp. $d\pi$)
with $t\circ\pi$ (resp. $d(t\circ\pi)$). 
Hence $\pi\in{\cal O}(\pi^{-1}({\cal U}))$
and $d\pi\in 
H^{0}(\pi^{-1}({\cal U}),\Omega^{1}_{X})$
in what follows.
\par
Let $h_{[\Gamma]}$ be a 
$C^{\infty}$ Hermitian metric on 
$[\Gamma]$ with
\begin{equation}
h_{[\Gamma]}(s_{\Gamma},s_{\Gamma})(w,t)
=
\begin{cases}
\begin{array}{lcr}
|\pi(w)-t|^{2}&\hbox{if}&(w,t)
\in
\pi^{-1}({\cal V})\times{\cal V},
\\
1&\hbox{if}&(w,t)
\in
(X\setminus\pi^{-1}({\cal U}))\times{\cal V}.
\end{array}
\end{cases}
\end{equation}
Let $h_{[\Gamma]^{-1}}$ be the metric
on $[\Gamma]^{-1}$ induced 
from $h_{[\Gamma]}$.
\par
\par
Let
$\|\cdot\|_{Q,\lambda(\xi)}$ 
be the Quillen metric on $\lambda(\xi)$ 
with respect to $g_{X/S}$, $h_{\xi}$.
Let
$\|\cdot\|_{Q,
\lambda([\Gamma]^{-1}\otimes p_{1}^{*}\xi)}$ 
(resp. $\|\cdot\|_{Q,\lambda(p_{1}^{*}\xi)}$)
be the Quillen metric on 
$\lambda([\Gamma]^{-1}\otimes p_{1}^{*}\xi)$ 
(resp. $\lambda(p_{1}^{*}\xi)$) 
with respect to 
$g_{X}$, $h_{[\Gamma]^{-1}}\otimes h_{\xi}$ 
(resp. $g_{X}$, $h_{\xi}$).
Let 
$\|\cdot\|_{Q,\lambda}$ 
be the Quillen metric on $\lambda$ 
defined as the tensor product of those on 
$\lambda([\Gamma]^{-1}\otimes p_{1}^{*}\xi)$,
$\lambda(p_{1}^{*}\xi)^{-1}$, 
$\lambda(\xi)$.
\par
For a complex manifold $Y$, 
$A^{p,q}(Y)$ denotes the vector space of
$C^{\infty}$ $(p,q)$-forms on $Y$.
We set
$\widetilde{A}(Y):=
\bigoplus_{p\geq0}A^{p,p}(Y)/
{\rm Im}\,\partial+{\rm Im}\,\bar{\partial}$.
\par
For a Hermitian vector bundle $(F,h_{F})$ 
over $Y$, 
$c_{i}(F,h_{F}),
{\rm Td}(F, h_{F}),
{\rm ch}(F,h_{F})\in
\bigoplus_{p\geq0}A^{p,p}(Y)$
denote the $i$-th Chern form,
the Todd form, and the Chern character form
of $(F,h_{F})$ with respect to 
the holomorphic Hermitian connection,
respectively. 
Let ${\rm R}(F)$ denote the ${\rm R}$-genus 
of Gillet-Soul\'e 
\cite[(0.4)]{BismutLebeau91},
\cite[p.\,160]{Soule92}.

\begin{theorem}
The following identity of functions on 
${\cal U}^{o}$ holds
$$
\log\|\sigma_{KM}\|^{2}_{Q,\lambda}
\equiv_{\cal B}
\left(\int_{E_{0}}
\widetilde{\mu}^{*}\left\{{\rm Td}(U)\,
\frac{{\rm Td}(H)-1}{c_{1}(H)}\right\}\,
q^{*}{\rm ch}(\xi)\right)\,\log|t|^{2}.
$$
\end{theorem}

\begin{pf}
We follow Bismut \cite[Sect.\,5]{Bismut97}.
(See also \cite[Th.\,6.3]{Yoshikawa04}.)
\newline{\bf (Step 1)}
Let $[X_{t}]$ be the holomorphic line bundle 
on $X$ associated to the divisor $X_{t}$. 
Then $[X_{t}]=[\Gamma]|_{X_{t}}$. 
We define the canonical section $s_{t}$ 
of $[X_{t}]$ by
$s_{t}:=s_{\Gamma}|_{X\times\{t\}}\in 
H^{0}(X,[X_{t}])$. 
Then ${\rm div}(s_{t})=X_{t}$. 
Let $i_{t}\colon X_{t}\hookrightarrow X$ 
be the embedding, 
and set $\xi_{t}:=\xi|_{X_{t}}$.
By (5.1), we get the exact sequence of 
coherent sheaves on $X$,
\begin{equation}
0
\longrightarrow
{\cal O}_{X}([X_{t}]^{-1}\otimes\xi)
@>\otimes s_{t}>>
{\cal O}_{X}(\xi)
\longrightarrow
(i_{t})_{*}{\cal O}_{X_{t}}(\xi)
\longrightarrow
0.
\end{equation}
\par
Let $\lambda([X_{t}]^{-1}\otimes\xi)$ and 
$\lambda(\xi_{t})$ be the determinants of 
the cohomology groups of 
$[X_{t}]^{-1}\otimes\xi$ and $\xi_{t}$, respectively. 
Then 
$\lambda_{t}
=
\lambda([X_{t}]^{-1}\otimes\xi)\otimes
\lambda(\xi)^{-1}\otimes\lambda(\xi_{t})$.
\par
Set 
$h_{[X_{t}]}=h_{[\Gamma]}|_{X\times\{t\}}$
for $t\in{\cal V}$. 
Then $h_{[X_{t}]}$ is a Hermitian metric 
on $[X_{t}]$. 
Let $h_{[X_{t}]}^{-1}$ be the Hermitian metric 
on $[X_{t}]^{-1}$ induced from $h_{[X_{t}]}$.
\par
Let $N_{t}=N_{X_{t}/X}$ 
(resp. $N_{t}^{*}=N^{*}_{X_{t}/X}$) 
be the normal (resp. conormal) bundle 
of $X_{t}$ in $X$. 
Then 
$d\pi|_{X_{t}}\in H^{0}(X_{t},N_{t}^{*})$ 
generates $N^{*}_{t}$ for $t\in{\cal U}^{o}$. 
Let $h_{N^{*}_{t}}$ be the Hermitian metric 
on $N^{*}_{t}$ defined by
\begin{equation}
h_{N^{*}_{t}}(d\pi|_{X_{t}},d\pi|_{X_{t}})=1.
\end{equation}
Let $h_{N_{t}}$ be the Hermitian metric 
on $N_{t}$ induced from $h_{N^{*}_{t}}$. 
Then we have the identity 
$c_{1}(N_{t},h_{N_{t}})=0$
for $t\in{\cal V}^{o}$.
\par
For 
$(w,t)\in\pi^{-1}({\cal U})\times{\cal U}$, 
set 
$$
\widetilde{s}_{\Gamma}(w,t)
=
\frac{s_{\Gamma}(w,t)}{\pi(w)-t}.
$$
Since $\pi(w)-t$ is a holomorphic function on 
$\pi^{-1}({\cal U})\times{\cal U}$ 
with divisor $\Gamma$, 
$\widetilde{s}_{\Gamma}$
is a nowhere vanishing holomorphic section of
$[\Gamma]|_{\pi^{-1}({\cal U})\times{\cal U}}$. 
Set 
$\widetilde{s}_{X_{t}}
=
\widetilde{s}_{\Gamma}|_{X_{t}\times\{t\}}
\in H^{0}(X_{t},[X_{t}]|_{X_{t}})$ 
and
$$
ds_{t}|_{X_{t}}
:=
d\pi\otimes\widetilde{s}_{X_{t}}
\in 
H^{0}(X_{t},N^{*}_{t}\otimes[X_{t}]|_{X_{t}}).
$$
By (5.2), (5.4), the isomorphism
$$
\otimes ds_{t}|_{X_{t}}
\colon
[X_{t}]^{-1}\otimes\xi|_{X_{t}}
\ni 
v\to ds_{t}|_{X_{t}}(v)
\in 
N_{t}^{*}\otimes\xi_{t}
$$ 
gives an isometry of 
holomorphic Hermitian vector bundles
$$
([X_{t}]^{-1}\otimes\xi,
h_{[X_{t}]^{-1}}\otimes h_{\xi})|_{X_{t}}
\cong
(N_{t}^{*}\otimes\xi_{t},
h_{N_{t}^{*}}\otimes h_{\xi}|_{X_{t}})
$$ 
for all $t\in{\cal V}^{o}$.  
Hence the metrics
$h_{[X_{t}]^{-1}}\otimes h_{\xi}$ 
and 
$h_{\xi}$ 
verify assumption (A) of 
Bismut \cite[Def.1.5]{Bismut90} 
with respect to 
$h_{N_{t}}$ and $h_{\xi}|_{X_{t}}$.
\newline{\bf (Step 2)}
Associated to the exact sequence of 
holomorphic vector bundles on $X_{t}$,
$$
{\cal E}_{t}
\colon
0
\longrightarrow 
TX_{t}
\longrightarrow 
TX|_{X_{t}}
\longrightarrow 
N_{t}
\longrightarrow
0,
$$
one can define the Bott-Chern class
$\widetilde{\rm Td}
({\cal E}_{t};g_{X_{t}},g_{X},h_{N_{t}})
\in\widetilde{A}(X_{t})$ 
by
\cite[I, f)]{BGS88}, 
\cite[I, Sect.\,1]{GilletSoule90},
\cite[Chap.\,IV, Sect.\,3]{Soule92}
such that
$$
dd^{c}\widetilde{\rm Td}
({\cal E}_{t};g_{X_{t}},g_{X},h_{N_{t}})
=
{\rm Td}(TX_{t},g_{X_{t}})\,
{\rm Td}(N_{t},h_{N_{t}})
-
{\rm Td}(TX,g_{X})|_{X_{t}}.
$$
Notice that our
$\widetilde{\rm Td}
({\cal E}_{t};\,g_{X_{t}},g_{X},h_{N_{t}})$
and Bismut-Lebeau's
$\widetilde{\rm Td}
(TX_{t},TX|_{X_{t}},h_{N_{t}})$ are related
as follows:
$$
\widetilde{\rm Td}
({\cal E}_{t};\,g_{X_{t}},g_{X},h_{N_{t}})
=
-\widetilde{\rm Td}
(TX_{t},TX|_{X_{t}},h_{N_{t}}).
$$
\par  
Let $Z$ be a general fiber of 
$\pi\colon X\to S$.
By applying the embedding formula of 
Bismut-Lebeau \cite[Th.\,0.1]{BismutLebeau91} 
(see also \cite[Th.\,5.6]{Bismut97}) to 
the embedding
$i_{t}\colon X_{t}\hookrightarrow X$ and 
to the exact sequence (5.3), we get 
for all $t\in{\cal V}^{o}$:
\begin{equation}
\begin{aligned}
\log\|\sigma_{KM}(t)\|^{2}_{Q,\lambda}
&
=
\int_{X\times\{t\}}
-\frac{{\rm Td}(TX,g_{X})\,
{\rm ch}(\xi,h_{\xi})}
{{\rm Td}([\Gamma],h_{[\Gamma]})}
\log h_{[\Gamma]}(s_{\Gamma},s_{\Gamma})
|_{X\times\{t\}}
\\
&
\quad
-\int_{X_{t}}\frac{\widetilde{\rm Td}
({\cal E}_{t};\,g_{X_{t}},g_{X},h_{N_{t}})\,
{\rm ch}(\xi,h_{\xi})}
{{\rm Td}(N_{t},h_{N_{t}})}
\\
&
\quad
-
\int_{X}{\rm Td}(TX)\,R(TX)\,{\rm ch}(\xi)
+
\int_{Z}{\rm Td}(TZ)\,
R(TZ)\,{\rm ch}(\xi|_{Z}).
\end{aligned}
\end{equation}
Here we used the explicit formula 
for the Bott-Chern current
\cite[Rem.\,3.5, especially (3.23),
Th.\,3.15, Th.\,3.17]{BGS90}
to get the first term of 
the right hand side of (5.5).
Notice that the dual of 
our $\lambda(\xi)$ was defined as
$\lambda(\xi)$ in \cite{BismutLebeau91}.
\par
By Theorem 9.1 below,
the first term of the right hand side of (5.5)
lies in ${\cal B}({\cal U})$.
Substituting $c_{1}(N_{t},h_{N_{t}})=0$ 
into (5.5), 
we get
\begin{equation}
\log\|\sigma_{KM}(t)\|^{2}_{Q,\lambda}
\equiv_{\cal B}
\int_{X_{t}}-\widetilde{\rm Td}
({\cal E}_{t};\,g_{X_{t}},g_{X},h_{N_{t}})\,
{\rm ch}(\xi,h_{\xi}).
\end{equation}
\newline{\bf (Step 3)}
Let $g_{N_{t}}$ be the Hermitian metric 
on $N_{t}$ induced from $g_{X}$ 
by the $C^{\infty}$ isomorphism
$N_{t}\cong(TX_{t})^{\perp}$. 
Let 
$\widetilde{\rm Td}
(N_{t};\,h_{N_{t}},g_{N_{t}})
\in\widetilde{A}(X_{t})$
be the Bott-Chern class 
\cite[I, e)]{BGS88}, 
\cite[Sect.\,1.2.4]{GilletSoule90},
\cite[Chap.\,IV, Sect.\,3]{Soule92}
such that
$$
dd^{c}
\widetilde{\rm Td}(N_{t};h_{N_{t}},g_{N_{t}})
=
{\rm Td}(N_{t},h_{N_{t}})
-{\rm Td}(N_{t},g_{N_{t}}).
$$
By \cite[I, Prop.\,1.3.2 and Prop.\,1.3.4]
{GilletSoule90}
(see also Lemma 5.3 below),
\begin{equation}
\widetilde{\rm Td}
({\cal E}_{t};\,g_{X_{t}},g_{X},h_{N_{t}})
=
\widetilde{\rm Td}
({\cal E}_{t};\,g_{X_{t}},g_{X},g_{N_{t}})
+
{\rm Td}(TX_{t},g_{X_{t}})\,
\widetilde{\rm Td}
(N_{t};\,h_{N_{t}},g_{N_{t}}).
\end{equation}
Since 
$c_{1}(N_{t},h_{N_{t}})=0$ 
and
$g_{N_{t}}=\|d\pi\|^{-2}\,h_{N_{t}}$, 
we deduce from
\cite[I, Prop.\,1.3.1 and (1.2.5.1)]
{GilletSoule90}
the identity
\begin{equation}
\begin{aligned}
\widetilde{\rm Td}(N_{t};\,h_{N_{t}},g_{N_{t}})
&
=
\frac{1-{\rm Td}(dd^{c}\log\|d\pi\|^{2})}
{dd^{c}\log\|d\pi\|^{2}}\,\log\|d\pi\|^{2}
\\
&
=
\left.\nu^{*}\left\{
\frac{1-{\rm Td}(-c_{1}(L,g_{L}))}
{-c_{1}(L,g_{L})}
\right\}\,
\log\|d\pi\|^{2}\right|_{X_{t}}.
\end{aligned}
\end{equation}
Substituting (5.8) and
$(TX_{t},g_{X_{t}})
=
\mu^{*}(U,g_{U})|_{X_{t}}$
into (5.7), we get
\begin{equation}
\begin{aligned}
\,&
\widetilde{\rm Td}
({\cal E}_{t};\,g_{X_{t}},g_{X},h_{N_{t}})
=
\\
&
\widetilde{\rm Td}
({\cal E}_{t};\,g_{X_{t}},g_{X},g_{N_{t}})
+
\left.\mu^{*}{\rm Td}(U,g_{U})\,
\nu^{*}\left\{
\frac{1-{\rm Td}(-c_{1}(L,g_{L}))}
{-c_{1}(L,g_{L})}
\right\}\,\log\|d\pi\|^{2}\right|_{X_{t}}.
\end{aligned}
\end{equation}
Since
$$
{\cal E}_{t}
=
\mu^{*}{\cal S}^{\lor}|_{X_{t}},
\quad
g_{X_{t}}=\mu^{*}g_{U}|_{X_{t}},
\quad
g_{X}
=
\mu^{*}(\varPi^{\lor})^{*}g_{X}|_{X_{t}},
\quad
g_{N_{t}}=\mu^{*}g_{H}|_{X_{t}},
$$ 
we deduce from
\cite[I, Th.\,1.2.2 (ii)]{GilletSoule90}
that
\begin{equation}
\widetilde{\rm Td}
({\cal E}_{t};\,g_{X_{t}},g_{X},g_{N_{t}})
=
\mu^{*}\widetilde{\rm Td}
({\cal S}^{\lor};\,
g_{U},(\varPi^{\lor})^{*}g_{X},g_{H})
|_{X_{t}}.
\end{equation}
Comparing (5.9) and (5.10), we get
\begin{equation}
\begin{aligned}
\widetilde{\rm Td}
({\cal E}_{t};\,g_{X_{t}},g_{X},h_{N_{t}})
&
=
\mu^{*}\widetilde{\rm Td}
({\cal S}^{\lor};\,
g_{U},(\varPi^{\lor})^{*}g_{X},g_{H})
|_{X_{t}}
\\
&
\quad
+
\mu^{*}{\rm Td}(U,g_{U})\,
\nu^{*}\left\{
\frac{1-{\rm Td}(-c_{1}(L,g_{L}))}
{-c_{1}(L,g_{L})}
\right\}\,
\log\|d\pi\|^{2}|_{X_{t}}.
\end{aligned}
\end{equation}
Substituting (5.11) into (5.6), we get
\begin{equation}
\begin{aligned}
\,&
\log\|\sigma_{KM}\|^{2}_{Q,\lambda}
\\
&
\equiv_{\cal B}
-\pi_{*}\left[\mu^{*}\widetilde{\rm Td}
({\cal S}^{\lor};\,
g_{U},(\varPi^{\lor})^{*}g_{X},g_{H})\,
{\rm ch}(\xi,h_{\xi})\right]^{(0,0)}
\\
&
\quad-
\pi_{*}\left[\mu^{*}{\rm Td}(U,g_{U})\,
\nu^{*}\left\{
\frac{1-{\rm Td}(-c_{1}(L,g_{L}))}
{-c_{1}(L,g_{L})}
\right\}\,{\rm ch}(\xi,h_{\xi})\,
\log\|d\pi\|^{2}
\right]^{(0,0)}
\\
&
\equiv_{\cal B}
-\widetilde{\pi}_{*}\left[
\widetilde{\mu}^{*}\widetilde{\rm Td}
({\cal S}^{\lor};\,
g_{U},(\varPi^{\lor})^{*}g_{X},g_{H})\,
q^{*}{\rm ch}(\xi,h_{\xi})\right]^{(0,0)}
\\
&
\quad+
\widetilde{\pi}_{*}\left[
\widetilde{\mu}^{*}{\rm Td}(U,g_{U})\,
\widetilde{\nu}^{*}\left\{
\frac{{\rm Td}(-c_{1}(L,g_{L}))-1}
{-c_{1}(L,g_{L})}
\right\}\,
q^{*}{\rm ch}(\xi,h_{\xi})\,
(q^{*}\log\|d\pi\|^{2})\right]^{(0,0)}.
\end{aligned}
\end{equation}
\par
Recall that
for a $C^{\infty}$ differential form 
$\varphi$ on $\widetilde{X}$, 
one has
$\widetilde{\pi}_{*}(\varphi)^{(0,0)}\in
{\cal B}({\cal U})$ by Barlet 
\cite[Th.\,4bis]{Barlet82}.
Since $q^{*}{\rm ch}(\xi,h_{\xi})$ and
$$
\widetilde{\mu}^{*}\widetilde{\rm Td}
({\cal S}^{\lor};\,
g_{U},(\varPi^{\lor})^{*}g_{X},g_{H}),
\quad
\widetilde{\mu}^{*}{\rm Td}(U,g_{U}),
\quad 
\widetilde{\nu}^{*}\{
\frac{{\rm Td}(-c_{1}(L,g_{L}))-1}
{-c_{1}(L,g_{L})}\}
$$ 
are $C^{\infty}$ differential forms 
on $\widetilde{X}$, 
we deduce from (5.12), 
\cite[Th.\,4bis]{Barlet82}, 
and Corollary 4.6 that
\begin{equation}
\log\|\sigma_{KM}\|^{2}_{Q,\lambda}
\equiv_{\cal B}
\left(\int_{E_{0}}
\widetilde{\mu}^{*}\left\{{\rm Td}(U)\,
\frac{{\rm Td}(H)-1}{c_{1}(H)}\right\}\,
q^{*}{\rm ch}(\xi)\right)\,\log|t|^{2}.
\end{equation}
Here we used the identity
$c_{1}(H)
=
-c_{1}(L)+(\varPi^{\lor})^{*}\pi^{*}c_{1}(S)$
in $H^{2}({\Bbb P}(TX)^{\lor},{\Bbb Z})$ and
the triviality of the line bundle
$\widetilde{\mu}^{*}(\varPi^{\lor})^{*}
\pi^{*}(TS)|_{\widetilde{\pi}^{-1}({\cal U})}$ 
to get (5.13). 
This completes the proof of Theorem 5.1.
\end{pf}

For simplicity, we set
$\overline{L}:=(L,g_{L})$,
$\overline{U}:=(U,g_{U})$,
$\overline{\xi}:=(\xi,h_{\xi})$
in what follows.
\par
Let
$\widetilde{\rm Td}({\cal S}^{\lor};
\,g_{U},(\varPi^{\lor})^{*}g_{X},g_{H})$ 
be the Bott-Chern secondary class associated with 
the Todd genus and the exact sequence of 
holomorphic vector bundles
$$
{\cal S}^{\lor}
\colon
0
\to 
U
\to
(\varPi^{\lor})^{*}TX
\to 
H
\to 
0
$$
equipped with the Hermitian metrics
$g_{U}$, $(\varPi^{\lor})^{*}g_{X}$, $g_{H}$, 
such that
$$
dd^{c}\widetilde{\rm Td}
({\cal S}^{\lor};\,
g_{U},(\varPi^{\lor})^{*}g_{X},g_{H})
=
{\rm Td}(U,g_{U})\,{\rm Td}(H,g_{H})
-
(\varPi^{\lor})^{*}{\rm Td}(TX,g_{X}).
$$
Recall that $Z$ is a general fiber of 
$\pi\colon X\to S$. 

\begin{theorem}
The following identity holds
$$
\begin{aligned}
\,&
\lim_{t\to0}
\left[
\log\|\sigma_{KM}(t)\|^{2}_{Q,\lambda}
-
\left(
\int_{E_{0}}\widetilde{\mu}^{*}
\left\{{\rm Td}(U)\,
\frac{{\rm Td}(H)-1}{c_{1}(H)}\right\}\,
q^{*}{\rm ch}(\xi)
\right)
\log\|{\bf s}_{0}(t)\|_{0}^{2}\right]
=
\\
&
-\int_{X\times\{0\}}
\frac{{\rm Td}(TX,g_{X})\,
{\rm ch}(\overline{\xi})}
{{\rm Td}([\Gamma],h_{[\Gamma]})}
\log \|s_{\Gamma}\|^{2}|_{X\times\{0\}}
\\
&
-\int_{\widetilde{X}_{0}}
\widetilde{\mu}^{*}\widetilde{\rm Td}
({\cal S}^{\lor};\,
g_{U},(\varPi^{\lor})^{*}g_{X},g_{H})\,
q^{*}{\rm ch}(\overline{\xi})
\\
&
+\int_{\widetilde{X}}
(q^{*}\log\|d\pi\|^{2})\,
\widetilde{\pi}^{*}c_{1}([0],\|\cdot\|_{0})
\left[
\widetilde{\mu}^{*}{\rm Td}(\overline{U})\,
\widetilde{\nu}^{*}\left\{
\frac{{\rm Td}(-c_{1}(\overline{L}))-1}
{-c_{1}(\overline{L})}
\right\}\,q^{*}{\rm ch}(\overline{\xi})
\right]
\\
&
-
\int_{\widetilde{X}}
(\widetilde{\pi}^{*}
\log\|{\bf s}_{0}\|_{0}^{2})\,
\widetilde{\nu}^{*}c_{1}(\overline{L})\,
\left[
\widetilde{\mu}^{*}{\rm Td}(\overline{U})\,
\widetilde{\nu}^{*}\left\{
\frac{{\rm Td}(-c_{1}(\overline{L}))-1}
{-c_{1}(\overline{L})}
\right\}\,q^{*}{\rm ch}(\overline{\xi})
\right]
\\
&
-\int_{X}{\rm Td}(TX)\,{\rm R}(TX)\,
{\rm ch}(\xi)
+\int_{Z}{\rm Td}(TZ)\,{\rm R}(TZ)\,
{\rm ch}(\xi|_{Z}).
\end{aligned}
$$
\end{theorem}

\begin{pf}
Define topological constants $C_{0}$ and 
$C_{1}$ by
$$
C_{0}
:=
\int_{E_{0}}\widetilde{\mu}^{*}
\left\{{\rm Td}(U)\,
\frac{{\rm Td}(H)-1}{c_{1}(H)}\right\}\,
q^{*}{\rm ch}(\xi),
$$
$$
C_{1}
:=
-\int_{X}{\rm Td}(TX)\,{\rm R}(TX)\,
{\rm ch}(\xi)
+\int_{Z}{\rm Td}(TZ)\,{\rm R}(TZ)\,
{\rm ch}(\xi|_{Z}).
$$
Substituting (5.11) and 
$c_{1}(N_{t},h_{N_{t}})=0$ into (5.5), 
we get for $t\in{\cal U}^{o}$
\begin{equation}
\begin{aligned}
\log\|\sigma_{KM}(t)\|^{2}_{Q,\lambda}
&
=
-\int_{X\times\{t\}}
\frac{{\rm Td}(TX,g_{X})\,
{\rm ch}(\overline{\xi})}
{{\rm Td}([\Gamma],h_{[\Gamma]})}
\log \|s_{\Gamma}\|^{2}
|_{X\times\{t\}}
\\
&
\quad
-\int_{X_{t}}
\mu^{*}\widetilde{\rm Td}
({\cal S}^{\lor};\,
g_{U},(\varPi^{\lor})^{*}g_{X},g_{H})
|_{X_{t}}\,
{\rm ch}(\overline{\xi})
\\
&
\quad
-
\int_{X_{t}}
\mu^{*}{\rm Td}(\overline{U})\,
\nu^{*}\left\{
\frac{1-{\rm Td}(-c_{1}(\overline{L}))}
{-c_{1}(\overline{L})}
\right\}\,
{\rm ch}(\overline{\xi})
\,\log\|d\pi\|^{2}
+C_{1}
\\
&
=
-\int_{X\times\{t\}}
\frac{{\rm Td}(TX,g_{X})\,
{\rm ch}(\overline{\xi})}
{{\rm Td}([\Gamma],h_{[\Gamma]})}
\log \|s_{\Gamma}\|^{2}
|_{X\times\{t\}}
\\
&
\quad
-\int_{\widetilde{X}_{t}}
\widetilde{\mu}^{*}\widetilde{\rm Td}
({\cal S}^{\lor};\,
g_{U},(\varPi^{\lor})^{*}g_{X},g_{H})|_{X_{t}}
\,q^{*}{\rm ch}(\overline{\xi})
\\
&
\quad
+
\int_{\widetilde{X}_{t}}
\widetilde{\mu}^{*}{\rm Td}(\overline{U})\,
\widetilde{\nu}^{*}\left\{
\frac{{\rm Td}(-c_{1}(\overline{L}))-1}
{-c_{1}(\overline{L})}
\right\}\,
q^{*}{\rm ch}(\overline{\xi})\,
q^{*}(\log\|d\pi\|^{2})
+C_{1},
\end{aligned}
\end{equation}
which yields that
\begin{equation}
\begin{aligned}
\,&
\log\|\sigma_{KM}(t)\|^{2}_{Q,\lambda}
-C_{0}\,\log\|{\bf s}_{0}(t)\|_{0}^{2}
=
\\
&
-\int_{X\times\{t\}}
\frac{{\rm Td}(TX,g_{X})\,
{\rm ch}(\overline{\xi})}
{{\rm Td}([\Gamma],h_{[\Gamma]})}
\log \|s_{\Gamma}\|^{2}
-
\int_{\widetilde{X}_{t}}
\widetilde{\mu}^{*}\widetilde{\rm Td}
({\cal S}^{\lor};\,
g_{U},(\varPi^{\lor})^{*}g_{X},g_{H})
\,q^{*}{\rm ch}(\overline{\xi})
\\
&
+
\int_{\widetilde{X}_{t}}
\left[
\widetilde{\mu}^{*}{\rm Td}(\overline{U})\,
\widetilde{\nu}^{*}\left\{
\frac{{\rm Td}(-c_{1}(\overline{L}))-1}
{-c_{1}(\overline{L})}
\right\}\,
q^{*}{\rm ch}(\overline{\xi})
\right]
\,q^{*}(\log\|d\pi\|^{2})
-C_{0}\,\log\|{\bf s}_{0}(t)\|_{0}^{2}
\\
&
+C_{1}.
\end{aligned}
\end{equation}
By Corollary 4.5, 
\begin{equation}
\begin{aligned}
\,&
\int_{\widetilde{X}_{t}}
\left[
\widetilde{\mu}^{*}{\rm Td}(\overline{U})\,
\widetilde{\nu}^{*}\left\{
\frac{{\rm Td}(-c_{1}(\overline{L}))-1}
{-c_{1}(\overline{L})}
\right\}\,q^{*}{\rm ch}(\overline{\xi})
\right]
\,q^{*}(\log\|d\pi\|^{2})
-C_{0}\,\log\|{\bf s}_{0}(t)\|_{0}^{2}
\\
&
=
\int_{\widetilde{X}}
(q^{*}\log\|d\pi\|^{2})\,
\widetilde{\pi}^{*}c_{1}([0],\|\cdot\|_{0})
\left[
\widetilde{\mu}^{*}{\rm Td}(\overline{U})\,
\widetilde{\nu}^{*}\left\{
\frac{{\rm Td}(-c_{1}(\overline{L}))-1}
{-c_{1}(\overline{L})}
\right\}\,q^{*}{\rm ch}(\overline{\xi})
\right]
\\
&
\quad
-
\int_{\widetilde{X}}
(\widetilde{\pi}^{*}
\log\|{\bf s}_{0}\|_{0}^{2})\,
\widetilde{\nu}^{*}c_{1}(\overline{L})\,
\left[
\widetilde{\mu}^{*}{\rm Td}(\overline{U})\,
\widetilde{\nu}^{*}\left\{
\frac{{\rm Td}(-c_{1}(\overline{L}))-1}
{-c_{1}(\overline{L})}
\right\}\,q^{*}{\rm ch}(\overline{\xi})
\right]+o(1).
\end{aligned}
\end{equation}
From (5.15) and (5.16), we get
\begin{equation}
\begin{aligned}
\,&
\lim_{t\to0}
\left[
\log\|\sigma_{KM}(t)\|^{2}_{Q,\lambda}
-C_{0}\,\log\|{\bf s}_{0}(t)\|_{0}^{2}
\right]
=
\\
&
-\int_{X\times\{0\}}
\frac{{\rm Td}(TX,g_{X})\,
{\rm ch}(\overline{\xi})}
{{\rm Td}([\Gamma],h_{[\Gamma]})}
\log \|s_{\Gamma}\|^{2}
|_{X\times\{0\}}
\\
&
-\int_{\widetilde{X}_{0}}
\widetilde{\mu}^{*}\widetilde{\rm Td}
({\cal S}^{\lor};\,
g_{U},(\varPi^{\lor})^{*}g_{X},g_{H})
\,q^{*}{\rm ch}(\overline{\xi})
\\
&
+\int_{\widetilde{X}}
(q^{*}\log\|d\pi\|^{2})\,
\widetilde{\pi}^{*}c_{1}([0],\|\cdot\|_{0})
\left[
\widetilde{\mu}^{*}{\rm Td}(\overline{U})\,
\widetilde{\nu}^{*}\left\{
\frac{{\rm Td}(-c_{1}(\overline{L}))-1}
{-c_{1}(\overline{L})}
\right\}\,q^{*}{\rm ch}(\overline{\xi})
\right]
\\
&
-
\int_{\widetilde{X}}
(\widetilde{\pi}^{*}
\log\|{\bf s}_{0}\|_{0}^{2})\,
\widetilde{\nu}^{*}c_{1}(\overline{L})\,
\left[
\widetilde{\mu}^{*}{\rm Td}(\overline{U})\,
\widetilde{\nu}^{*}\left\{
\frac{{\rm Td}(-c_{1}(\overline{L}))-1}
{-c_{1}(\overline{L})}
\right\}\,q^{*}{\rm ch}(\overline{\xi})
\right]
+C_{1}.
\end{aligned}
\end{equation}
This completes the proof of Theorem 5.2.
\end{pf}

\begin{lemma}
Let 
$
{\cal E}\colon
0
\longrightarrow
E'
\longrightarrow
E
\longrightarrow
E''
\longrightarrow
0
$
be an exact sequence of holomorphic
vector bundles over a complex manifold
$Y$. 
Let $h'$ and $h$ be Hermitian metrics 
on $E'$ and $E$, respectively.
Let $h''$ and $g''$ be Hermitian metrics
on $E''$. Then
$$
\widetilde{\rm Td}({\cal E};\,h',h,h'')
-
\widetilde{\rm Td}({\cal E};\,h',h,g'')
=
{\rm Td}(E',h')\,
\widetilde{\rm Td}(E'';\,h'',g'').
$$
\end{lemma}

\begin{pf}
Setting
$\overline{L}_{1}=({\cal E},h',h,h'')$,
$\overline{L}_{2}=({\cal E},h',h,g'')$,
$\overline{L}_{3}=0$ 
in \cite[I, Prop.\,1.3.4]{GilletSoule90},
we get
$$
\widetilde{\rm Td}({\cal E};\,h',h,h'')
-
\widetilde{\rm Td}({\cal E};\,h',h,g'')
=
\widetilde{\rm Td}(E'\oplus E'';\,
h'\oplus h'',h'\oplus g'').
$$
Since
$\widetilde{\rm Td}(E'\oplus E'';\,
h'\oplus h'',h'\oplus g'')
=
{\rm Td}(E',h')\,
\widetilde{\rm Td}(E'';\,h'',g'')$
by \cite[I, Prop.\,1.3.2]{GilletSoule90},
we get the result.
\end{pf}


\section
{\bf The divergent term and the constant term}
\par

Let $\alpha$ be a nowhere vanishing
holomorphic section of
$\lambda([\Gamma]^{-1}\otimes p_{1}^{*}\xi)^{-1}
\otimes\lambda(p_{1}^{*}\xi)$ defined on ${\cal U}$.

\begin{theorem}
Let $\sigma$ be a nowhere vanishing 
holomorphic section of $\lambda(\xi)$ 
defined on $\cal U$. 
Then
$$
\log\|\sigma\|^{2}_{Q,\lambda(\xi)}
\equiv_{\cal B}
\left(\int_{E_{0}}
\widetilde{\mu}^{*}\left\{{\rm Td}(U)\,
\frac{{\rm Td}(H)-1}{c_{1}(H)}\right\}\,
q^{*}{\rm ch}(\xi)\right)\,\log|t|^{2}.
$$
\end{theorem}

\begin{pf}
There exists a nowhere vanishing holomorphic
function $f(t)$ on $\cal U$ such that
$$
\sigma(t)
=
f(t)\,\sigma_{KM}(t)\otimes\alpha(t).
$$
Since $\log|f(t)|^{2}$ and
$\log\|\alpha\|^{2}_{Q,
\lambda([\Gamma]^{-1}\otimes p_{1}^{*}\xi)^{-1}
\otimes\lambda(p_{1}^{*}\xi)}$ are $C^{\infty}$ 
functions on $\cal U$, we deduce from 
Theorem 5.1 that
$$
\begin{aligned}
\log\|\sigma(t)\|^{2}_{Q,\lambda(\xi)}
&=
\log|f(t)|^{2}+
\log\|\sigma_{KM}(t)\|^{2}_{Q,\lambda}+
\log\|\alpha(t)\|^{2}_{Q,
\lambda([\Gamma]^{-1}\otimes p_{1}^{*}\xi)^{-1}
\otimes\lambda(p_{1}^{*}\xi)}
\\
&\equiv_{\cal B}
\left(\int_{E_{0}}
\widetilde{\mu}^{*}\left\{{\rm Td}(U)\,
\frac{{\rm Td}(H)-1}{c_{1}(H)}\right\}\,
q^{*}{\rm ch}(\xi)\right)\,\log|t|^{2}.
\end{aligned}
$$
This completes the proof of Theorem 6.1.
\end{pf}

\begin{theorem}
The following identity holds:
$$
\begin{aligned}
\,&
\lim_{t\to0}
\left[
\log\|
\sigma_{KM}\otimes\alpha
\|^{2}_{Q,\lambda(\xi)}(t)-
\left(
\int_{E_{0}}\widetilde{\mu}^{*}
\left\{{\rm Td}(U)\,
\frac{{\rm Td}(H)-1}{c_{1}(H)}\right\}\,
q^{*}{\rm ch}(\xi)
\right)
\log\|{\bf s}_{0}(t)\|_{0}^{2}\right]
\\
&=
\log\|\alpha(0)\|_{Q}^{2}
-\int_{X\times\{0\}}
\frac{{\rm Td}(TX,g_{X})\,
{\rm ch}(\overline{\xi})}
{{\rm Td}([\Gamma],h_{[\Gamma]})}
\log \|s_{\Gamma}\|^{2}|_{X\times\{0\}}
\\
&
\quad
-\int_{\widetilde{X}_{0}}
\widetilde{\mu}^{*}\widetilde{\rm Td}
({\cal S}^{\lor};\,
g_{U},(\varPi^{\lor})^{*}g_{X},g_{H})\,
q^{*}{\rm ch}(\overline{\xi})
\\
&
\quad
+\int_{\widetilde{X}}
(q^{*}\log\|d\pi\|^{2})\,
\widetilde{\pi}^{*}c_{1}([0],\|\cdot\|_{0})
\left[
\widetilde{\mu}^{*}{\rm Td}(\overline{U})\,
\widetilde{\nu}^{*}\left\{
\frac{{\rm Td}(-c_{1}(\overline{L}))-1}
{-c_{1}(\overline{L})}
\right\}\,q^{*}{\rm ch}(\overline{\xi})
\right]
\\
&
\quad
-
\int_{\widetilde{X}}
(\pi^{*}\log\|{\bf s}_{0}\|_{0}^{2})\,
\widetilde{\nu}^{*}c_{1}(\overline{L})\,
\left[
\widetilde{\mu}^{*}{\rm Td}(\overline{U})\,
\widetilde{\nu}^{*}\left\{
\frac{{\rm Td}(-c_{1}(\overline{L}))-1}
{-c_{1}(\overline{L})}
\right\}\,q^{*}{\rm ch}(\overline{\xi})
\right]
\\
&
\quad
-\int_{X}{\rm Td}(TX)\,{\rm R}(TX)\,
{\rm ch}(\xi)
+\int_{Z}{\rm Td}(TZ)\,{\rm R}(TZ)\,
{\rm ch}(\xi|_{Z}).
\end{aligned}
$$
\end{theorem}

\begin{pf}
Since
$$
\log\|\sigma_{KM}\otimes\alpha
\|_{Q,\lambda(\xi)}^{2}
=
\log\|\sigma_{KM}\|_{Q,\lambda}^{2}
+
\log\|\alpha\|^{2}_{Q,
\lambda([\Gamma]^{-1}\otimes p_{1}^{*}\xi)^{-1}
\otimes\lambda(p_{1}^{*}\xi)},
$$
the result follows from Theorem 5.2.
\end{pf}


\section
{\bf Critical points defined by 
a quadric polynomial of rank $2$}
\par
In this section, we assume that for every
$x\in\Sigma_{\pi}\cap X_{0}$, 
there exists a system of coordinates 
$(z_{0},\ldots,z_{n})$ centered at $x$ 
such that
$$
\pi(z)=z_{0}z_{1}.
$$
Hence $\Sigma_{\pi}\subset X$ is 
a complex submanifold of codimension $2$
defined locally by the equation $z_{0}=z_{1}=0$. 
Let 
$N_{\Sigma_{\pi}/X}$ be the normal bundle 
of $\Sigma_{\pi}$ in $X$. 
In \cite[Def.\,5.1, Prop.\,5.2]{Bismut97},
Bismut introduced the additive genus 
${\rm E}(\cdot)$ 
associated with the generating function
$$
{\rm E}(x)
:=
\frac{{\rm Td}(x)\,{\rm Td}(-x)}{2x}
\left(\frac{{\rm Td}^{-1}(x)-1}{x}
-\frac{{\rm Td}^{-1}(-x)-1}{-x}\right),
$$
where
${\rm Td}^{-1}(x):=(1-e^{-x})/x$.
\par
The following result was proved
by Bismut \cite[Th.\,5.9]{Bismut97}.

\begin{theorem}
The following equation of functions 
on ${\cal U}^{o}$ holds:
$$
\log
\|\sigma(t)\|^{2}_{\lambda(\xi),Q}
\equiv_{\cal B}
\frac{1}{2}\,
\left(\int_{\Sigma_{\pi}\cap X_{0}}
-{\rm Td}(T\Sigma_{\pi})\,
{\rm E}(N_{\Sigma_{\pi}/X})
\,{\rm ch}(\xi)\right)\,\log|t|^{2}.
$$
\end{theorem}

\begin{remark}
As mentioned before, the dual of
our $\lambda(\xi)$ was defined as
$\lambda(\xi)$ in \cite[Th.\,5.9]{Bismut97},
which explains the difference of 
the sign of the coefficient of $\log|t|^{2}$ 
in Theorem 7.1 with that of 
\cite[Th.\,5.9]{Bismut97}.
\end{remark}

\begin{pf}
Let $q\colon\widetilde{X}\to X$ be 
the blowing-up along $\Sigma_{\pi}$ 
with exceptional divisor
$$
E={\Bbb P}(N_{\Sigma_{\pi}/X}).
$$
Then $\widetilde{\nu}=\nu\circ q$ extends 
to a holomorphic map from $\widetilde{X}$ to 
${\Bbb P}(\Omega^{1}_{X})$. 
\par
Since the Hessian of $\pi$ is 
a non-degenerate symmetric bilinear form 
on $N_{\Sigma_{\pi}/X}$, 
we have
$N_{\Sigma_{\pi}/X}\cong 
N_{\Sigma_{\pi}/X}^{*}$. 
Under the identification
${\Bbb P}(N_{\Sigma_{\pi}/X})
=
{\Bbb P}(N_{\Sigma_{\pi}/X}^{*})$ 
induced from the Hessian of $\pi$, 
$\widetilde{\nu}$ is identified with 
the natural inclusion
${\Bbb P}(N_{\Sigma_{\pi}/X}^{*})
\hookrightarrow{\Bbb P}(\Omega^{1}_{X}
|_{\Sigma_{\pi}})$, 
which yields that
\begin{equation}
\widetilde{\nu}^{*}L|_{E}
=
{\cal O}_{{\Bbb P}(N_{\Sigma_{\pi}/X}^{*})}(-1),
\qquad
\widetilde{\mu}^{*}H|_{E}
=
{\cal O}_{{\Bbb P}(N_{\Sigma_{\pi}/X})}(1).
\end{equation}
Set 
$F
:=
{\cal O}_{{\Bbb P}(N_{\Sigma_{\pi}/X})}(1)$.
\par
By the exact sequence ${\cal S}^{\lor}$, 
we get
\begin{equation}
{\rm Td}(U)=
\frac{{\rm Td}((\varPi^{\lor})^{*}TX)}
{{\rm Td}(H)}.
\end{equation}
Since $\varPi^{\lor}\circ\widetilde{\mu}=q$,
we deduce from the exact sequence of vector 
bundles on $\Sigma_{\pi}$
$$
0\longrightarrow
T\Sigma_{\pi}\longrightarrow
TX|_{\Sigma_{\pi}}\longrightarrow
N_{\Sigma_{\pi}/X}\longrightarrow0
$$
the identity
\begin{equation}
\widetilde{\mu}^{*}
{\rm Td}((\varPi^{\lor})^{*}TX)
|_{E}
=
q^{*}\left\{{\rm Td}(T\Sigma_{\pi})\,
{\rm Td}(N_{\Sigma_{\pi}/X})\right\}.
\end{equation}
Substituting (7.3) into (7.2), 
we get
\begin{equation}
\widetilde{\mu}^{*}{\rm Td}(U)|_{E}
=
\frac{q^{*}\left\{{\rm Td}(T\Sigma_{\pi})\,
{\rm Td}(N_{\Sigma_{\pi}/X})\right\}}
{\widetilde{\mu}^{*}{\rm Td}(H)|_{E}}
=
\frac{q^{*}\left\{{\rm Td}(T\Sigma_{\pi})\,
{\rm Td}(N_{\Sigma_{\pi}/X})\right\}}
{{\rm Td}(F)},
\end{equation}
where we used (7.1) to get the second equality.
\par
Let $p_{*}$ be the integration 
along the fibers of the projection
$p\colon{\Bbb P}(N_{\Sigma_{\pi}/X})\to
\Sigma_{\pi}$.
Since $q|_{E}=p$, we deduce from
(7.1), (7.4) and the projection formula that
\begin{equation}
\begin{aligned}
\,&
\int_{E\cap X_{0}}
\widetilde{\mu}^{*}\left\{{\rm Td}(U)\,
\frac{{\rm Td}(H)-1}{c_{1}(H)}\right\}\,
q^{*}{\rm ch}(\xi)
\\
&
=
\int_{\Sigma_{\pi}\cap X_{0}}
{\rm Td}(T\Sigma_{\pi})\,
{\rm Td}(N_{\Sigma_{\pi}/X})\,{\rm ch}(\xi)\,
p_{*}\left\{\frac{1}{{\rm Td}(F)}\cdot
\frac{{\rm Td}(F)-1}{c_{1}(F)}\right\}
\\
&
=
\int_{\Sigma_{\pi}\cap X_{0}}
{\rm Td}(T\Sigma_{\pi})\,
{\rm Td}(N_{\Sigma_{\pi}/X})\,{\rm ch}(\xi)\,
p_{*}\left\{
\frac{1-{\rm Td}^{-1}(F)}{c_{1}(F)}
\right\}.
\end{aligned}
\end{equation}
Since 
$N_{\Sigma_{\pi}/X}\cong
N_{\Sigma_{\pi}/X}^{*}$, 
we have
$$
c_{1}(N_{\Sigma_{\pi}/X})=0,
$$
which, together with
${\rm rk}(N_{\Sigma_{\pi}/X})=2$, 
yields that
$$
0
=
c_{1}(F)^{2}-
p^{*}c_{1}(N_{\Sigma_{\pi}/X})\,c_{1}(F)+
p^{*}c_{2}(N_{\Sigma_{\pi}/X})
=
c_{1}(F)^{2}+p^{*}c_{2}(N_{\Sigma_{\pi}/X}).
$$
Since $p_{*}c_{1}(F)=1$, this implies that 
for $m\geq0$
\begin{equation}
p_{*}c_{1}(F)^{m}=
\begin{cases}
\begin{array}{lr}
(-1)^{k}\,c_{2}(N_{\Sigma_{\pi}/X})^{k}
&(m=2k+1)
\\
0&(m=2k).
\end{array}
\end{cases}
\end{equation}
\par
For a formal power series
$f(x)=\sum_{j=0}^{\infty}a_{j}\,x^{j}\in{\Bbb C}[[x]]$, 
set
$$
f_{-}(x):=\frac{f(x)-f(-x)}{2x}
\in{\Bbb C}[[x]].
$$ 
By (7.6), we get
$$
p_{*}f(c_{1}(F))
=
\sum_{k}a_{2k+1}\,p_{*}c_{1}(F)^{2k+1}
=
\sum_{k}(-1)^{k}a_{2k+1}\,
c_{2}(N_{\Sigma_{\pi}/X})^{k}.
$$
Let $f_{-}(N_{\Sigma_{\pi}/X})$ be 
the additive genus associated with 
$f_{-}(x)\in{\Bbb C}[[x]]$.
Let $x_{1}$, $x_{2}$ be the Chern roots 
of $N_{\Sigma_{\pi}/X}$. 
Since 
$c_{1}(N_{\Sigma_{\pi}/X})=x_{1}+x_{2}=0$,
we get
$$
\begin{aligned}
f_{-}(N_{\Sigma_{\pi}/X})
&
=
\frac{f(x_{1})-f(-x_{1})}{2x_{1}}
+
\frac{f(x_{2})-f(-x_{2})}{2x_{2}}
\\
&
=
\sum_{k=0}^{\infty}a_{2k+1}\,
(x_{1}^{2k}+x_{2}^{2k})
\\
&
=
2\sum_{k=0}^{\infty}a_{2k+1}(-x_{1}x_{2})^{k}
\\
&
=
2\sum_{k=0}^{\infty}(-1)^{k}a_{2k+1}\,
c_{2}(N_{\Sigma_{\pi}/X})^{k}
=
2\,p_{*}f(c_{1}(F)).
\end{aligned}
$$
Setting
$f(x)=({\rm Td}^{-1}(x)-1)/x$, 
we get
\begin{equation}
\begin{aligned}
{\rm E}(N_{\Sigma_{\pi}/X})
&
=
{\rm Td}(x_{1}){\rm Td}(x_{2})\,
\left\{
\frac{f(x_{1})-f(-x_{1})}{2x_{1}}
+
\frac{f(x_{2})-f(-x_{2})}{2x_{2}}
\right\}
\\
&
=
2\,{\rm Td}(N_{\Sigma_{\pi}/X})\,
p_{*}f(c_{1}(F))
\\
&
=
-2\,{\rm Td}(N_{\Sigma_{\pi}/X})\,
p_{*}\left(
\frac{1-{\rm Td}^{-1}(F)}{c_{1}(F)}
\right).
\end{aligned}
\end{equation}
By comparing (7.5) and (7.7), 
the desired formula follows from Theorem 6.1.
\end{pf}


\section
{\bf Isolated critical points}
\par
In this section, we assume that 
${\rm Sing}(X_{0})=\Sigma_{\pi}\cap X_{0}$ 
consists of isolated points. 
Since $\Sigma_{\pi}$ is discrete, 
we may identify
${\Bbb P}(\Omega^{1}_{X})$ and ${\Bbb P}(TX)$ 
with the trivial projective-space bundle 
on a neighborhood of $\Sigma_{\pi}\cap X_{0}$
by fixing a system of coordinates near 
$\Sigma_{\pi}\cap X_{0}$. 
Under this trivialization, 
we consider the Gauss maps $\nu$ and $\mu$ 
only on a small neighborhood of 
$\Sigma_{\pi}\cap X_{0}$. 
Then we have the following expression 
on a neighborhood of each 
$p\in\Sigma_{\pi}\cap X_{0}$:
$$
\mu(z)
=
\nu(z)
=
\left(
\frac{\partial\pi}{\partial z_{0}}(z)
:\cdots:
\frac{\partial\pi}{\partial z_{n}}(z)
\right).
$$
\par
For a formal power series 
$f(x)\in{\Bbb C}[[x]]$,
let $f(x)|_{x^{m}}$ denote
the coefficient of $x^{m}$.
Let $\mu(\pi,p)\in{\Bbb N}$ be 
the Milnor number of 
the isolated critical point $p$ of $\pi$.
The following result was proved
by the author \cite[Main Th.]{Yoshikawa98}.

\begin{theorem}
The following identity of functions 
on ${\cal U}^{o}$ holds:
$$
\log
\|\sigma\|^{2}_{\lambda(\xi),Q}
\equiv_{\cal B}
\frac{(-1)^{n}}{(n+2)!}\,{\rm rk}(\xi)
\left(\sum_{p\in{\rm Sing}(X_{0})}
\mu(\pi,p)\right)\,\log|t|^{2}.
$$
\end{theorem}

\begin{pf}
In Theorem 6.1, 
we can identify $U$ (resp. $L$)
with the universal hyperplane bundle 
(resp. tautological line bundle) 
on ${\Bbb P}^{n}$.
Then $H=L^{-1}$. 
Set $x:=c_{1}(H)$. 
Hence
$\int_{{\Bbb P}^{n}}x^{n}=1$.
From the exact sequence 
$0\to U\to{\Bbb C}^{n+1}\to H\to 0$, 
we get
$$
{\rm Td}(U)={\rm Td}^{-1}(x)=
\frac{1-e^{-x}}{x}.
$$
By substituting this and the equation
$q^{*}{\rm ch}(\xi)|_{E\cap\widetilde{X}_{0}}
={\rm rk}(\xi)$
into the formula of Theorem 6.1, 
we get
\begin{equation}
\begin{aligned}
\,&
\int_{E_{0}}
\widetilde{\mu}^{*}{\rm Td}(U)\,
\widetilde{\nu}^{*}
\left\{
\frac{{\rm Td}(c_{1}(H))-1}{c_{1}(H)}
\right\}\,
q^{*}{\rm ch}(\xi)
\\
&
=
\left.
\frac{1}{{\rm Td}(x)}\cdot
\frac{{\rm Td}(x)-1}{x}
\right|_{x^{n}}
\cdot
{\rm rk}(\xi)\int_{E_{0}}
\widetilde{\mu}^{*}c_{1}(H)^{n}
\\
&
=
\left.
\left\{\frac{1}{x}-\frac{1-e^{-x}}{x^{2}}
\right\}\right|_{x^{n}}
\cdot{\rm rk}(\xi)
\int_{E_{0}}\widetilde{\mu}^{*}c_{1}(H)^{n}
\\
&
=
\frac{(-1)^{n}}{(n+2)!}\,{\rm rk}(\xi)\,
\int_{E_{0}}\widetilde{\mu}^{*}c_{1}(H)^{n}.
\end{aligned}
\end{equation}
Since
$$
\begin{aligned}
\widetilde{\pi}_{*}\left\{
\widetilde{\mu}^{*}(-c_{1}(L,g_{L}))^{n}\,
q^{*}(\log\|d\pi\|^{2})\right\}
&
=
\pi_{*}\left\{
q^{*}(\log\|d\pi\|^{2})\,
(dd^{c}\log\|d\pi\|^{2})^{n}\right\}
\\
&
=
\sum_{p\in{\rm Sing}(X_{0})}
\mu(\pi,p)\,\log|t|^{2}
+O(1)
\end{aligned}
$$
by \cite[Th.\,4.1]{Yoshikawa98}, 
we get
\begin{equation}
\int_{E_{0}}\widetilde{\mu}^{*}c_{1}(H)^{n}
=
\sum_{p\in{\rm Sing}(X_{0})}\mu(\pi,p)
\end{equation}
by Corollary 4.6.
The result follows from Theorem 6.1 
and (8.1), (8.2).
\end{pf}


\section
{\bf Some results on asymptotic expansion}
\par

Let ${\cal A}_{\Bbb C}$ 
(resp. ${\cal C}_{\Bbb C}$) be 
the sheaf of germs of $C^{\infty}$ 
(resp. $C^{0}$) functions on $\Bbb C$. 
The stalk of ${\cal A}_{\Bbb C}$ 
(resp. ${\cal C}_{\Bbb C}$)
at the origin is denoted by 
${\cal A}_{0}$ (resp. ${\cal C}_{0}$).
We define
$$
{\cal B}_{0}
:=
{\cal A}_{0}
\oplus
\bigoplus_{r\in{\Bbb Q}\cap(0,1]}
\bigoplus_{k=0}^{n}
|t|^{2r}(\log|t|)^{k}\cdot{\cal A}_{0}
\subset
{\cal C}_{0}.
$$
In this section, we prove the following

\begin{theorem}
Let $\varOmega\subset{\Bbb C}^{n}$ be 
a relatively compact domain.
Let $F(z)$ be a holomorphic function 
on $\varOmega$ with critical locus
$\Sigma_{F}:=\{z\in\varOmega;\,dF(z)=0\}$.
Let $\chi(z)$ be a $C^{\infty}$ $(n,n)$-form
with compact support in $\varOmega$.
Define a germ $\psi\in{\cal C}_{0}$ by
$$
\psi(t)
:=
\int_{\varOmega}\log|F(z)-t|^{2}\,\chi(z).
$$
If $\Sigma_{F}\subset F^{-1}(0)$,
then $\psi(t)\in{\cal B}_{0}$.
\end{theorem}

The continuity of similar integrals
was studied by Bost-Gillet-Soul\'e
\cite[Sect.\,1.5]{BostGilletSoule94}
in relation with the regularity
of the star products of Green currents.
\par
For the proof of Theorem 9.1, we prove
some intermediary results.

\begin{lemma}
Let $\Phi$ be a $C^{\infty}$ $(n,n)$-form 
with compact support in $\varOmega$. 
Let $F_{*}(\Phi)$ be the locally integrable 
$(1,1)$-form on $\Bbb C$ defined as 
the integration of $\Phi$ 
along the fibers of 
$F\colon\varOmega\to{\Bbb C}$. 
If $\Sigma_{F}\subset F^{-1}(0)$,
then there exists a germ 
$A(t)\in{\cal B}_{0}$
such that
$$
F_{*}(\Phi)(t)
=
A(t)\,\frac{dt\wedge d\bar{t}}{|t|^{2}},
\qquad 
A(0)=0
$$
near $0\in{\Bbb C}$.
\end{lemma}

\begin{pf}
By Hironaka, 
there exists a proper holomorphic modification 
$\varpi\colon
\widetilde{\varOmega}\to\varOmega$ 
such that
\newline{(i)}
$\varpi\colon 
\widetilde{\varOmega}\setminus
\varpi^{-1}(\Sigma_{F})
\to 
\varOmega\setminus\Sigma_{F}$ 
is an isomorphism;
\newline{(ii)}
$(F\circ\varpi)^{-1}(\Sigma_{F})$ is 
a normal crossing divisor 
of $\widetilde{\varOmega}$.
\par
Set $\widetilde{F}:=F\circ\varpi$.
For any $z\in F^{-1}(0)$, 
there exist a system of coordinates 
$(U,(w_{1},\ldots,w_{n}))$ and
integers $k_{1},\ldots,k_{l}\geq1$,
$l\leq n$, 
such that 
$\widetilde{F}(w)
=
w_{1}^{k_{1}}\cdots w_{l}^{k_{l}}$. 
Define a holomorphic $(n-1)$-form on $U$
by
$$
\tau
:=
\frac{1}{l}\sum_{i=1}^{l}
\frac{1}{k_{i}}(-1)^{i-1}
w_{i}\,dw_{1}\wedge\cdots\wedge dw_{i-1}
\wedge dw_{i+1}\wedge\cdots\wedge dw_{n}.
$$
\par
Let $\varrho_{U}$ be a $C^{\infty}$ function 
with compact supported in $U$.
Since $\varpi^{*}\Phi$ is 
a $C^{\infty}$ $(n,n)$-form 
on $\widetilde{\varOmega}$, 
there exists $h(w)\in C^{\infty}_{0}(U)$ 
such that
$$
\varrho_{U}\varpi^{*}\Phi
=
h(w)\,dw_{1}\wedge\cdots\wedge dw_{n}\wedge 
d\bar{w}_{1}\wedge\cdots\wedge d\bar{w}_{n}.
$$
We define a germ $B(t)\in{\cal C}_{0}$ by
$$
B(t)
:=
\int_{\widetilde{F}^{-1}(t)\cap U}
h(w)\,\tau\wedge\bar{\tau}.
$$
Then $B(t)\in{\cal B}_{0}$ by
\cite[p.166, Th.\,4bis]{Barlet82}.
Since 
$$
\widetilde{F}^{*}\left(\frac{dt}{t}\right)
\wedge\tau
=
dw_{1}\wedge\cdots\wedge dw_{n},
$$ 
we get by the projection formula
\begin{equation}
\begin{aligned}
\widetilde{F}_{*}
(\varrho_{U}\,\varpi^{*}\Phi)(t)
&
=
\widetilde{F}_{*}
(h(w)\,dw_{1}\wedge\cdots\wedge dw_{n}\wedge 
d\bar{w}_{1}\wedge\cdots\wedge d\bar{w}_{n})(t)
\\
&
=
\frac{dt\wedge d\bar{t}}{|t|^{2}}\,
\widetilde{F}_{*}\left(
h(w)\,\tau\wedge\bar{\tau}\right)
=
B(t)\,\frac{dt\wedge d\bar{t}}{|t|^{2}}.
\end{aligned}
\end{equation}
\par
For an $\epsilon>0$ small enough, 
set 
$\varDelta(\epsilon)
:=
\{t\in{\Bbb C};\,|t|<\epsilon\}$. 
Since
$$
\left|
\int_{\varDelta(\epsilon)}
\widetilde{F}_{*}(\varrho_{U}\,
\varpi^{*}\Phi)
\right|
=
\left|
\int_{\widetilde{F}^{-1}(\varDelta(\epsilon))}
\varrho_{U}\,\varpi^{*}\Phi
\right|
<
\infty,
$$
the $(1,1)$-form 
$B(t)\,dt\wedge d\bar{t}/|t|^{2}$
is locally integrable near the origin.
Hence $B(0)=0$.
\par
Let $\{U_{\beta}\}_{\beta\in B}$ be
a locally finite open covering of 
$\widetilde{\varOmega}$
and
let $\{\varrho_{\beta}\}_{\beta\in B}$ be
a partition of unity subject to
$\{U_{\beta}\}_{\beta\in B}$. 
By (9.1), there exists 
$B_{\beta}(t)\in{\cal B}_{0}$ 
for each $\beta\in B$
such that
$$
\widetilde{F}_{*}(\varrho_{\beta}\,
\varpi^{*}(\Phi))
=
B_{\beta}(t)\,
\frac{dt\wedge d\bar{t}}{|t|^{2}},
\qquad
B_{\beta}(0)=0.
$$
There exist finitely many $\beta\in B$
with $B_{\beta}(t)\not=0$ by the compactness
of the support of $\varpi^{*}\Phi$.
Since
$$
F_{*}(\Phi)
=
\sum_{\beta\in B}\widetilde{F}_{*}
(\varrho_{\beta}\,\varpi^{*}\Phi)
=
(\sum_{\beta\in B}B_{\beta}(t))\,
\frac{dt\wedge d\bar{t}}{|t|^{2}},
$$
we get
$A(t)=\sum_{\beta\in B}B_{\beta}(t)
\in{\cal B}_{0}$ and $A(0)=0$.
\end{pf}

We regard $\varOmega$ 
as a domain in $({\Bbb P}^{1})^{n}$.
Hence $\chi$ is a $C^{\infty}$ $(n,n)$-form
on $({\Bbb P}^{1})^{n}$.
Let $z=(z_{1},\ldots,z_{n})$ be 
the inhomogeneous coordinates of
$({\Bbb P}^{1})^{n}$.
For $1\leq i\leq n$, 
set
$$
\omega_{i}
:=
\frac{\sqrt{-1}\,dz_{i}\wedge d\bar{z}_{i}}
{2\pi(1+|z_{i}|^{2})^{2}}.
$$

\begin{lemma}
Assume that
$F(z)=z_{1}^{\nu_{1}}\cdots z_{n}^{\nu_{n}}$,
$\nu_{1},\ldots,\nu_{n}\geq0$
and set
$$
\alpha
:=
\int_{({\Bbb P}^{1})^{n}}
\chi(z).
$$
Then there exists $\eta(t)\in{\cal B}_{0}$ 
such that
$$
\psi(t)
=
\alpha\int_{({\Bbb P}^{1})^{n}}
\log|
z_{1}^{\nu_{1}}\cdots z_{n}^{\nu_{n}}-t
|^{2}\,
\omega_{1}\wedge\cdots\wedge\omega_{n}
+
\eta(t).
$$
\end{lemma}

\begin{pf}
Let 
$((\zeta_{1}:\xi_{1}),\ldots,
(\zeta_{n}:\xi_{n}))$ be 
the homogeneous coordinates of 
$({\Bbb P}^{1})^{n}$
such that
$z_{i}=\zeta_{i}/\xi_{i}$.
For $t\in{\Bbb C}$, set
$$
Y_{t}
:=
\{
((\zeta_{1}:\xi_{1}),\ldots,
(\zeta_{n}:\xi_{n}))\in
({\Bbb P}^{1})^{n};\,
\zeta_{1}^{\nu_{1}}\cdots\zeta_{n}^{\nu_{n}}
-
t\,\xi_{1}^{\nu_{1}}\cdots\xi_{n}^{\nu_{n}}
=0
\},
$$
$$
D
:=
\{
((\zeta_{1}:\xi_{1}),\ldots,
(\zeta_{n}:\xi_{n}))\in
({\Bbb P}^{1})^{n};\,
\xi_{1}^{\nu_{1}}\cdots\xi_{n}^{\nu_{n}}
=0
\}.
$$
Since 
\begin{equation}
z_{1}^{\nu_{1}}\cdots z_{n}^{\nu_{n}}-t
=
\frac{\zeta_{1}^{\nu_{1}}\cdots\zeta_{n}^{\nu_{n}}
-t\,\xi_{1}^{\nu_{1}}\cdots\xi_{n}^{\nu_{n}}}
{\xi_{1}^{\nu_{1}}\cdots\xi_{n}^{\nu_{n}}},
\end{equation}
we get the following equation of currents
on $({\Bbb P}^{1})^{n}$ 
by the Poincar\'e-Lelong formula:
\begin{equation}
dd^{c}\log|
z_{1}^{\nu_{1}}\cdots z_{n}^{\nu_{n}}-t
|^{2}
=
\delta_{Y_{t}}-\delta_{D}.
\end{equation}
\par
Since $\chi(z)$ is cohomologous to
$\alpha\,
\omega_{1}\wedge\cdots\wedge\omega_{n}$,
there exists a $C^{\infty}$ $(n-1,n-1)$-form
$\gamma$ on $({\Bbb P}^{1})^{n}$ 
by the $dd^{c}$-Poincar\'e lemma, 
such that
$$
\chi(z)
-
\alpha\,\omega_{1}\wedge\cdots\wedge\omega_{n}
=
dd^{c}\gamma.
$$
Hence we get by (9.3)
\begin{equation}
\begin{aligned}
\psi(t)
&
=
\alpha\int_{({\Bbb P}^{1})^{n}}
\log|
z_{1}^{\nu_{1}}\cdots z_{n}^{\nu_{n}}-t
|^{2}\,
\omega_{1}\wedge\cdots\wedge\omega_{n}
+
\int_{({\Bbb P}^{1})^{n}}
\log|
z_{1}^{\nu_{1}}\cdots z_{n}^{\nu_{n}}-t
|^{2}\,
dd^{c}\gamma
\\
&
=
\alpha\int_{({\Bbb P}^{1})^{n}}
\log|
z_{1}^{\nu_{1}}\cdots z_{n}^{\nu_{n}}-t
|^{2}\,
\omega_{1}\wedge\cdots\wedge\omega_{n}
+
\int_{({\Bbb P}^{1})^{n}}
dd^{c}(\log|
z_{1}^{\nu_{1}}\cdots z_{n}^{\nu_{n}}-t
|^{2})\wedge\gamma
\\
&
=
\alpha\int_{({\Bbb P}^{1})^{n}}
\log|
z_{1}^{\nu_{1}}\cdots z_{n}^{\nu_{n}}-t
|^{2}\,
\omega_{1}\wedge\cdots\wedge\omega_{n}
+
\int_{Y_{t}}\gamma
-
\int_{D}\gamma.
\end{aligned}
\end{equation}
\par
For $t\in{\Bbb C}$, set
\begin{equation}
\eta(t)
:=
\int_{Y_{t}}\gamma-\int_{D}\gamma.
\end{equation}
Define a divisor of
$({\Bbb P}^{1})^{n}\times{\Bbb C}$ by
$$
Y
:=
\{
((\zeta_{1}:\xi_{1}),\ldots,
(\zeta_{n}:\xi_{n}),t)\in
({\Bbb P}^{1})^{n}\times{\Bbb C};\,
\zeta_{1}^{\nu_{1}}\cdots\zeta_{n}^{\nu_{n}}
-
t\xi_{1}^{\nu_{1}}\cdots\xi_{n}^{\nu_{n}}
=0
\}.
$$
Let 
${\rm pr}_{1}\colon
({\Bbb P}^{1})^{n}\times{\Bbb C}
\to({\Bbb P}^{1})^{n}$
and
${\rm pr}_{2}\colon
({\Bbb P}^{1})^{n}\times{\Bbb C}
\to{\Bbb C}$ 
be the projections. 
Then $Y_{t}=Y\cap{\rm pr}_{2}^{-1}(t)$.
Let $P\colon\widetilde{Y}\to Y$
be the resolution of the singularities
of $Y$. Then
${\rm pr}_{2}|_{Y}\circ P$
is a proper holomorphic function on 
the complex manifold $\widetilde{Y}$. 
Since
$P^{*}({\rm pr}_{1})^{*}\gamma$
is a $C^{\infty}$ $(n-1,n-1)$-form on
$\widetilde{Y}$, we get
\begin{equation}
\eta(t)
=
\int_{({\rm pr}_{2}|_{Y}\circ P)^{-1}(t)}
P^{*}({\rm pr}_{1})^{*}\gamma
-
\int_{D}\gamma
\in
{\cal B}_{0}
\end{equation}
by \cite[Th.\,4bis]{Barlet82}.
The result follows from (9.4), (9.5), (9.6).
\end{pf}

Define a germ $f\in{\cal C}_{0}$ by
$$
f(t)
:=
\int_{({\Bbb P}^{1})^{n}}
\log|
z_{1}^{\nu_{1}}\cdots z_{n}^{\nu_{n}}-t
|^{2}\,
\omega_{1}\wedge\cdots\wedge\omega_{n}.
$$

\begin{lemma}
There exists a germ $g(t)\in{\cal B}_{0}$ 
such that
$$
dd^{c}f(t)
=
\frac{\sqrt{-1}}{4\pi}\,
g(t)\,
\frac{dt\wedge d\bar{t}}{|t|^{2}},
\qquad
g(0)=0.
$$
\end{lemma}

\begin{pf}
We keep the notation 
in the proof of Lemma 9.3.
Since the assertion is obvious 
when $\nu_{1}=\cdots=\nu_{n}=0$,
we assume that $\nu_{i}>0$ for some $i$.
Since 
$z_{1}^{\nu_{1}}\cdots z_{n}^{\nu_{n}}-t$
is a meromorphic function 
on $({\Bbb P}^{1})^{n}\times{\Bbb C}$,
we deduce from (9.2) and 
the Poincar\'e-Lelong formula 
the following equation of currents 
on $({\Bbb P}^{1})^{n}\times{\Bbb C}$:
\begin{equation}
dd^{c}\log|
z_{1}^{\nu_{1}}\cdots z_{n}^{\nu_{n}}-t
|^{2}
=
\delta_{Y}-\delta_{D\times{\Bbb C}}
=
\delta_{Y}-\delta_{({\rm pr}_{1})^{*}D}.
\end{equation}
Since
$$
f
=
({\rm pr}_{2})_{*}
\left\{
\log|
z_{1}^{\nu_{1}}\cdots z_{n}^{\nu_{n}}-t
|^{2}\,
({\rm pr}_{1})^{*}
(\omega_{1}\wedge\cdots\wedge\omega_{n})
\right\},
$$
we get on ${\Bbb C}\setminus\{0\}$
\begin{equation}
\begin{aligned}
dd^{c}f
&
=
({\rm pr}_{2})_{*}
\left\{
dd^{c}\log|
z_{1}^{\nu_{1}}\cdots z_{n}^{\nu_{n}}-t
|^{2}\wedge
({\rm pr}_{1})^{*}
(\omega_{1}\wedge\cdots\wedge\omega_{n})
\right\}
\\
&
=
({\rm pr}_{2})_{*}
\left\{
(\delta_{Y}-\delta_{({\rm pr}_{1})^{*}D})
\wedge({\rm pr}_{1})^{*}
(\omega_{1}\wedge\cdots\wedge\omega_{n})
\right\}
\\
&
=
({\rm pr}_{2})_{*}
\left\{
({\rm pr}_{1})^{*}
(\omega_{1}\wedge\cdots\wedge\omega_{n})|_{Y}
\right\}
-
({\rm pr}_{2})_{*}
\left\{
({\rm pr}_{1})^{*}
(\omega_{1}\wedge\cdots\wedge\omega_{n}
|_{D})
\right\}
\\
&
=
({\rm pr}_{2}|_{Y})_{*}
\left\{
({\rm pr}_{1})^{*}
(\omega_{1}\wedge\cdots\wedge\omega_{n})|_{Y}
\right\}
\\
&
=
({\rm pr}_{2}|_{Y}\circ P)_{*}
\left\{
P^{*}({\rm pr}_{1})^{*}
(\omega_{1}\wedge\cdots\wedge\omega_{n})
\right\},
\end{aligned}
\end{equation}
where the first equality follows from 
the commutativity
$dd^{c}({\rm pr}_{2})_{*}
=
({\rm pr}_{2})_{*}dd^{c}$, 
the second equality follows from (9.7),
and the fourth equality follows from 
the trivial identity
$\omega_{1}\wedge\cdots\wedge\omega_{n}
|_{D}=0$.
Since 
$P^{*}({\rm pr}_{1})^{*}
(\omega_{1}\wedge\cdots\wedge\omega_{n})$ 
is a $C^{\infty}$ $(n,n)$-form
on $\widetilde{Y}$
and since 
${\rm pr}_{2}|_{Y}\circ P\colon
\widetilde{Y}\to{\Bbb C}$ is
a proper holomorphic map,
the assertion follows from 
(9.8) and Lemma 9.2.
\end{pf}

\begin{lemma}
The germ $f(t)$ is $S^{1}$-invariant,
i.e., $f(t)=f(|t|)$.
\end{lemma}

\begin{pf}
Without loss of generality,
we may assume that $\nu_{n}>0$.
Since
$$
\int_{{\Bbb P}^{1}}
\log|Az_{n}^{\nu_{n}}+B|^{2}\,\omega_{n}
=
\log(|A|^{2/\nu_{n}}+|B|^{2/\nu_{n}})
$$
when $(A,B)\not=(0,0)$, 
we get by Fubini's theorem
\begin{equation}
\begin{aligned}
f(t)
&
=
\int_{({\Bbb P}^{1})^{n}}
\log|
z_{1}^{\nu_{1}}\cdots z_{n}^{\nu_{n}}-t
|^{2}\,
\omega_{1}\wedge\cdots\wedge\omega_{n}
\\
&
=
\int_{({\Bbb P}^{1})^{n-1}}
\left(
\int_{{\Bbb P}^{1}}
\log|
z_{1}^{\nu_{1}}\cdots z_{n}^{\nu_{n}}-t
|^{2}\,
\omega_{n}
\right)
\omega_{1}\wedge\cdots\wedge\omega_{n-1}
\\
&
=
\int_{({\Bbb P}^{1})^{n-1}}
\log
\left(
|z_{1}^{\nu_{1}}\cdots z_{n-1}^{\nu_{n-1}}
|^{2/\nu_{n}}
+|t|^{2/\nu_{n}}
\right)
\omega_{1}\wedge\cdots\wedge\omega_{n-1}.
\end{aligned}
\end{equation}
The assertion follows from (9.9).
\end{pf}

Let $(r,\theta)$ be the polar coordinates
of $\Bbb C$. Hence $t=r\,e^{i\theta}$.

\begin{lemma}
Let $\lambda(t)\in C^{\infty}(\varDelta^{*})$. 
Assume that $\lambda(t)$ is $S^{1}$-invariant, i.e.,
$\lambda(t)=\lambda(r)$.
If $r\,\partial_{r}\lambda(t)\in{\cal B}_{0}$,
then $\lambda(t)\in{\cal B}_{0}$.
\end{lemma}

\begin{pf}
By the definition of ${\cal B}_{0}$,
there exist a finite set
$A\subset{\Bbb Q}\cap(0,1]$
and germs
$\mu_{\alpha,k}(t)\in{\cal A}_{0}$,
$\alpha\in A$, $0\leq k\leq n$ such that
\begin{equation}
r\,\partial_{r}\lambda(r)
=
\sum_{\alpha\in A}\sum_{k=0}^{n}
r^{2\alpha}(\log r)^{k}\,\mu_{\alpha,k}(t).
\end{equation}
We may assume that
$\mu_{\alpha,k}(t)\in 
C^{\infty}(\varDelta(2\epsilon))$
for some $\epsilon>0$.
Since the left hand side of (9.10) is 
$S^{1}$-invariant, we may assume that
$\mu_{\alpha,k}(t)=\mu_{\alpha,k}(r)$
for all $\alpha$ and $k$ after replacing
$\mu_{\alpha,k}(t)$ by 
$\int_{0}^{2\pi}\mu_{\alpha,k}(e^{i\theta}t)\,
d\theta/2\pi$. 
By (9.10), we get
\begin{equation}
\lambda(\epsilon)-\lambda(r)
=
\sum_{\alpha\in A}\sum_{k=0}^{n}
\int_{r}^{\epsilon}
u^{2\alpha-1}(\log u)^{k}\,
\mu_{\alpha,k}(u)\,du.
\end{equation}
By (9.11), 
we see that $\lambda(t)\in{\cal C}_{0}$
by setting
$$
\lambda(0)
:=
\lambda(\epsilon)
-
\sum_{\alpha\in A}\sum_{k=0}^{n}
\int_{0}^{\epsilon}
u^{2\alpha-1}(\log u)^{k}\,
\mu_{\alpha,k}(u)\,du.
$$
Since $\lambda(t)\in{\cal C}_{0}$, 
we get by (9.11)
$$
\begin{aligned}
\lambda(r)
&
=
\lambda(0)
+
\sum_{\alpha\in A}\sum_{k=0}^{n}
\int_{0}^{r}u^{2\alpha-1}(\log u)^{k}
\mu_{\alpha,k}(u)\,du
\\
&
=
\lambda(0)
+
\sum_{\alpha\in A}\sum_{k=0}^{n}
r^{2\alpha}
\int_{0}^{1}v^{2\alpha-1}(\log r+\log v)^{k}
\mu_{\alpha,k}(vr)\,dv
\\
&
=
\lambda(0)
+
\sum_{\alpha\in A}\sum_{k=0}^{n}
\sum_{l=0}^{k}
\binom{k}{l}r^{2\alpha}(\log r)^{l}
\int_{0}^{1}v^{2\alpha-1}(\log v)^{k-l}
\mu_{\alpha,k}(vt)\,dv,
\end{aligned}
$$
which implies that
$\lambda(t)\in{\cal B}_{0}$.
\end{pf}

\begin{lemma}
If 
$F(z)=z_{1}^{\nu_{1}}\cdots z_{n}^{\nu_{n}}$,
$\nu_{1},\ldots,\nu_{n}\geq0$,
then $\psi(t)\in{\cal B}_{0}$.
\end{lemma}

\begin{pf}
By Lemma 9.3, it suffices to prove that
$f\in{\cal B}_{0}$.
Since $f(t)=f(r)$ by Lemma 9.5,
we deduce from Lemma 9.4 the equation
$$
\frac{1}{2\pi}
\partial_{t}\partial_{\bar{t}}f(t)
=
\frac{1}{4\pi}\{f''(r)+r^{-1}f'(r)\}
=
\frac{g(t)}{4\pi r^{2}}.
$$
Hence $g(t)$ is invariant under the rotation,
i.e., $g(t)=g(r)$, 
and the following equation holds
\begin{equation}
(r\,\partial_{r})^{2}f(r)=g(r).
\end{equation}
Since $g(t)\in{\cal B}_{0}$,
we deduce from Lemma 9.6 and (9.12) 
that
$r\,\partial_{r}f(r)\in{\cal B}_{0}$. 
By Lemma 9.6 again,
we get $f(t)\in{\cal B}_{0}$.
\end{pf}

{\bf Proof of Theorem 9.1 }
\newline
We keep the notation in the proof of
Lemma 9.2. There exists a system of
coordinate neighborhoods
$\{(U_{\beta},
w_{\beta}=(w_{1,\beta},\ldots,w_{n,\beta})
)\}_{\beta\in B}$ 
of $\widetilde{\varOmega}$
and integers
$k_{1,\beta},\ldots,k_{n,\beta}\geq0$ 
for each $\beta\in B$
such that 
$\widetilde{F}|_{U_{\beta}}(w_{\beta})
=
w_{1,\beta}^{k_{1,\beta}}\cdots 
w_{n,\beta}^{k_{n,\beta}}$.
Without loss of generality, we may assume that
the covering $\{U_{\beta}\}_{\beta\in B}$ 
of $\widetilde{\varOmega}$ is locally finite.
Let $\{\varrho_{\beta}\}_{\beta\in B}$ be
a partition of unity subject to the covering
$\{U_{\beta}\}_{\beta\in B}$.
Then
$\chi_{\beta}:=\varrho_{\beta}\,\varpi^{*}\chi$ 
is a $C^{\infty}$ $(n,n)$-form with compact support
in $U_{\beta}$. Since $\varpi^{*}\chi$ has
a compact support in $\widetilde{\Omega}$,
$\chi_{\beta}=0$ except finitely many $\beta\in B$.
By Lemma 9.7,
\begin{equation}
\psi_{\beta}(t)
:=
\int_{U_{\beta}}
\log|w_{1,\beta}^{k_{1,\beta}}\cdots
w_{n,\beta}^{k_{n,\beta}}-t|^{2}\,
\chi_{\beta}(w_{\beta})
\in{\cal B}_{0}.
\end{equation}
Since
$$
\psi(t)
=
\int_{\widetilde{\Omega}}\varpi^{*}
\log|F-t|^{2}\,\varpi^{*}\chi
=
\sum_{\beta\in B}\int_{U_{\beta}}
\log\left|
\widetilde{F}|_{U_{\beta}}(w_{\beta})-t\right|^{2}\,
\varrho_{\beta}\varpi^{*}\chi
=
\sum_{\beta\in B}\psi_{\beta}(t),
$$
we get $\psi(t)\in{\cal B}_{0}$ by (9.13).
This completes the proof of Theorem 9.1.
$\Box$


\end{document}